\documentclass[reqno, 11pt]{amsart}
\usepackage[utf8]{inputenc}
\usepackage[T1]{fontenc}
\usepackage{lmodern}
\usepackage{epsfig,amsmath, amsthm ,amsfonts,latexsym,stmaryrd, amssymb}
\usepackage[abs]{overpic}
\usepackage[usenames,dvipsnames]{xcolor}
\usepackage[linktocpage=true]{hyperref}
\usepackage{caption}
\usepackage{subcaption}
\usepackage{dsfont}
\usepackage{mathtools}
\usepackage{tikz-cd}
\usepackage{simpler-wick}
\usetikzlibrary{shapes,arrows}
\usepackage{comment}
\usepackage{enumitem}
\usepackage{marginnote}
\usepackage[all,cmtip]{xy}
\usepackage[backend=biber, maxnames=4, style=alphabetic]{biblatex}
\usepackage{tcolorbox}
\usepackage{graphicx}

\DeclareUnicodeCharacter{0306}{~}

\hypersetup{
  colorlinks,
  citecolor=teal,
  linkcolor=purple,
  urlcolor=teal}
\theoremstyle{plain}

\newtheorem{theorem}{Theorem}

\newtheorem{dummy}{anything}[section]

\theoremstyle{definition}

\newtheorem{example}[dummy]{Example}

\newtheorem{remark}[dummy]{Remark}

\theoremstyle{remark}
\newcommand{\R}{\mathbb{R}}

\def\R{\mathbb{R}}

\DeclareFontFamily{U}{mathx}{}
\DeclareFontShape{U}{mathx}{m}{n}{<-> mathx10}{}
\DeclareSymbolFont{mathx}{U}{mathx}{m}{n}
\DeclareMathAccent{\widecheck}{0}{mathx}{"71}
\newenvironment{nouppercase}{%
  \renewcommand{\uppercasenonmath}[1]{}}{}

\begin{document}

\title{\Large Symplectic geometry and space mission design \\ \vspace{0.5cm} \small On the Jupiter--Europa and Saturn--Enceladus systems}

\author{Cengiz Aydin}

\address[C.\ Aydin]{Universit\"at Heidelberg \\ Mathematisches Institut \\ Heidelberg \\ Germany}

\email{\href{mailto:cengiz.aydin@hotmail.de}{cengiz.aydin@hotmail.de}}

\author{Urs Frauenfelder}

\address[U.\ Frauenfelder]{Institut f\"ur Mathematik\\ Augsburg Universität\\ Augsburg \\ Germany} 

\email{\href{mailto:urs.frauenfelder@math.uni-augsburg.de}{urs.frauenfelder@math.uni-augsburg.de}}

\author{Otto van Koert}

\address[O.\ van Koert]{Dept.\ of Mathematical Sciences\\ Seoul National University \\ South Korea}

\email{\href{mailto:okoert@snu.ac.kr}{okoert@snu.ac.kr}}

\author{Dayung Koh} 

\address[D.\ Koh]{Pasadena\\ California}

\email{\href{mailto:dayung.koh@gmail.com}{dayung.koh@gmail.com}}

\author{Agustin Moreno}

\address[A.\ Moreno]{Universit\"at Heidelberg \\ Mathematisches Institut \\ Heidelberg \\ Germany}

\email{\href{mailto:agustin.moreno2191@gmail.com}{agustin.moreno2191@gmail.com}}

%\date{\today}

%\subjclass[2020]{Primary 57R17; Secondary 53C15, 32Q65, 37D05}

%57R17 Symplectic and contact topology
%53C15 General geometric structures on manifolds (almost complex, almost product structures, etc.)
%32Q65 Pseudoholomorphic curves
%37D05 Dynamical systems with hyperbolic orbits and sets
%%%%%%%%%%%%%%%%%%%%%%%%%%%%%%%%%%%%%%%%%%%%%%%%%%%%%%%%%%%%%%%%%%%%%
%%%%%%%%%%%%%%%%%%%%%%%% BEGIN PAPER %%%%%%%%%%%%%%%%%%%%%%%%%%%%%%%%
%%%%%%%%%%%%%%%%%%%%%%%%%%%%%%%%%%%%%%%%%%%%%%%%%%%%%%%%%%%%%%%%%%%%%

\begin{abstract}
   Symplectic geometry is a modern mathematical framework in which to address problems in mechanics, such as gravitational attraction and planetary motion. Using methods from this field, the second and fifth authors have provided theoretical groundwork and tools aimed at analyzing periodic orbits, their stability and their bifurcations in families, for the purpose of space mission design \cite{FM}. The Broucke stability diagram \cite{Broucke} was refined, and the "Floer numbers" were introduced, as numbers which stay invariant before and after a bifurcation, and therefore serve as tests for the algorithms used, as well as being easy to implement. These tools were further employed for numerical studies \cite{FKM}. In this article, we will further illustrate these methods with numerical studies of families of orbits for the Jupiter--Europa and Saturn--Enceladus systems, with emphasis on planar-to-spatial bifurcations, from deformation of the families in Hill's lunar problem studied by the first author \cite{Cengiz}. We will also provide an algorithm for the numerical computation of Conley--Zehnder indices, which are instrumental in practice for determining which families of orbits connect to which. 
\end{abstract}

\begin{nouppercase}
\maketitle    
\end{nouppercase}

\tableofcontents

\section{Introduction}
Symplectic geometry is the branch of mathematics that studies the geometric properties of phase spaces, those spaces that describe the possible states of a classical physical system. It provides a proper framework to address problems in classical mechanics, e.g.,\ the gravitational problem of N bodies in three-dimensional space. In the last thirty years, a host of theoretical tools have been developed in the field, with \emph{Floer theory} as a notable example, whose emphasis is on the theoretical study of periodic orbits. In a more applied direction, periodic orbits are of interest for space mission design, as they model trajectories for spacecraft or satellites. Studying families of orbits aimed at placing a spacecraft around a target moon is relevant for space exploration, where optimizing over all possible trajectories is needed, in order to minimize fuel consumption, avoid collisions, and maximize safety. In this context, the influence on a satellite of a planet with an orbiting moon can be approximated by a three-body problem of \emph{restricted} type (i.e.,\ the mass of the satellite is considered negligible by comparison). This is a classical problem which has been central to the development of symplectic geometry, and therefore it is not unreasonable to expect the modern available tools to provide insights. The need of organizing all information pertaining to orbits leads to the realm of data analysis, for which computationally cheap methods are important. The direction we will pursue is then encapsulated in the following questions:

\medskip

\begin{tcolorbox} \textbf{Guiding questions}

\medskip

    \begin{itemize}
    \item \textbf{(Classification)} Can we tell when two orbits are \emph{qualitatively different}? 

\medskip
    
    \item\textbf{(Catalogue)} Can we resource-efficiently refine data bases of known orbits? 

    \medskip
    
    \item \textbf{(Symplectic geometry)} Can we use methods from symplectic geometry to guide/organize the numerical work?
    \end{itemize}
\end{tcolorbox}

Here, we say that two orbits are \emph{qualitatively different} if they cannot be joined by a regular family of orbits, i.e.,\ a family which does \emph{not} undergo bifurcation. The first two questions were addressed by the second and fifth authors in \cite{FM}, where the mathematical groundwork was developed, and obstructions to the existence of regular families were encoded in the topology of suitable quotients of the symplectic group. This method, whose main tool is the \emph{GIT sequence}, gives a refinement of the well-known Broucke stability diagram \cite{Broucke}. This method was further developed for the case of Hamiltonian systems of arbitrary degrees of freedom by the fifth author and Ruscelli in \cite{MR}. The second, fourth, and fifth authors used it in combination with numerical work, addressing the third question \cite{FKM}. In this article, we continue this line of research. As before, we have the following tools at our disposal.

%\newpage

\begin{tcolorbox} \textbf{Toolkit}

\medskip
    
\begin{enumerate}
    
    \item[(1)] \textbf{Floer numerical invariants:} Integers which stay invariant before and after a bifurcation, and so can help predict the existence of orbits, as well as being easy to implement. There is one invariant for arbitrary periodic orbits, and another for \emph{symmetric} periodic orbits \cite{FKM}. 

\medskip
    
    \item[(2)] \textbf{The B-signs} \cite{FM}: a $\pm$ sign associated to each elliptic or hyperbolic Floquet multiplier of an orbit\footnote{Recall that the Floquet multipliers of a closed orbit are by definition the non-trivial eigenvalues of the monodromy matrix.}, which helps predict bifurcations. This is generalization of the classical Moser--Krein signature \cite{Kre2,Kre3,Moser}, which originally applies only to elliptic Floquet multipliers, to also include the case of hyperbolic multipliers, whenever the corresponding orbit is \emph{symmetric}. 

\medskip
    
    \item[(3)] \textbf{Global topological methods:} the \emph{GIT-sequence} \cite{FM}, a sequence of spaces whose global topology encodes (and sometimes forces) bifurcations, and refines Broucke's stability diagram \cite{Broucke} by adding the $B$-signs.

\medskip
    
    \item[(4)]\textbf{Conley-Zehnder index:} \cite{CZ,RS} a winding number associated to each non-degenerate orbit, extracted from the topology of the symplectic group, which does not change unless a bifurcation occurs. It can be used to determine which families connect to which.
\end{enumerate}

\end{tcolorbox}

Combining the Floer numbers with the $B$-signs provides tools to decide whether to look for periodic orbits, and gives hints concerning where to actually look for them. In this paper, which is an extension of \cite{AFKM}, we apply these tools in numerical studies of families of periodic orbits in the Saturn--Enceladus and the Jupiter--Europa system, by deformation of families in the lunar problem studied by the first author \cite{Cengiz}. Our results illustrate the general principle that one may learn about a given system, by starting from known nearby systems, and then deforming. One of the highlights of our paper are bifurcation graphs relating various families (Figures \ref{fig:bifurc_diagram_new} and \ref{fig:SE_graph}), including a \emph{spatial} family connecting two \emph{planar} orbits, one retrograde, and the other, prograde. We further provide the documentation for a numerical implementation of the CZ-indices in Appendix~\ref{app:CZ_index_numerical}, made publicly available via GitHub. We expect this to be useful in practice for the early stages of space mission design, where mapping out large data bases of periodic orbits is important. Finally, in Appendix \ref{app:halo}, we apply our methods in order to analyze a family of Halo orbits in the Saturn--Enceladus system, which approaches the plumes at an altitude of 29 km, and therefore may be used for future missions. Our novel tools were instrumental for our results. 

\medskip

\textbf{Acknowledgments.} A.\ Moreno received support by the NSF under Grant No.\ DMS-1926686, and is currently supported by the Sonderforschungsbereich TRR 191 Symplectic Structures in Geometry, Algebra and Dynamics, funded by the DFG (Projektnummer 281071066 – TRR 191), and by the DFG under Germany's Excellence Strategy EXC 2181/1 - 390900948 (the Heidelberg STRUCTURES Excellence Cluster).\
O.\ van Koert was supported by National Research Foundation of Korea Grant NRF2023005562 funded by the Korean Government.\
C.\ Aydin acknowledges support by the Deutsche Forschungsgemeinschaft (DFG, German Research Foundation) – Project-ID 281071066 – TRR 191.

\section{Preliminaries}

In this section, we review the toolkit.\ But first, we set up some language and notation. We refer the reader to \cite{FKM} for details on the global topological methods. 

\subsection{Basic notions} \textbf{Mechanics/symplectic geometry.} Given a $2n$-dimensional phase-space $M$ with its symplectic form $\omega$, a Hamiltonian function $H:M\rightarrow \mathbb R$, with Hamiltonian flow $\phi^H_t: M\rightarrow M$ which preserves $\omega$ (i.e.,\ $(\phi^H_t)^*\omega=\omega$), and a periodic orbit $x$, the \emph{monodromy matrix} of $x$ is $M_x=D\phi^H_T$, where $T$ is the period of $x$. Then $M_x$ is a symplectic $2n\times 2n$-matrix; we denote by $Sp(2n)$ the space of such matrices (the \emph{symplectic group}).

Note that if $H$ is time-independent then $1$ appears twice as a \emph{trivial} eigenvalue of $M_x$. We can ignore these if we consider the \emph{reduced} monodromy matrix $M_x^{red}\in Sp(2n-2)$, obtained by fixing the energy and dropping the direction of the flow.
    
\begin{itemize}
    \item A \emph{Floquet multiplier} of $x$ is an eigenvalue of $M_x$, which is not one of the trivial eigenvalues (i.e.,\ an eigenvalue of $M_x^{red}$).
    \item An orbit is \emph{non-degenerate} if $1$ does not appear among its Floquet multipliers. 
    \item An orbit is \emph{stable} if all its Floquet multipliers are semi-simple and lie on the unit circle.
\end{itemize}

We will only consider the cases $n=2$ (planar problems) and $n=3$ (spatial problems). 

\medskip

\textbf{Symmetries.} An \emph{anti-symplectic involution} is a map $\rho:M\rightarrow M$ satisfying $\rho^2=id$ and $\rho^*\omega=-\omega$. Its \emph{fixed-point locus} is $fix(\rho)=\{x:\rho(x)=x\}$. An anti-symplectic involution $\rho$ is a \emph{symmetry} of the system if $H\circ \rho=H.$ A periodic orbit $x$ is \emph{symmetric} if $\rho(x(-t))=x(t)$ for all $t$. The \emph{symmetric points} of the symmetric orbit $x$ are the two intersection points of $x$ with $fix(\rho)$. The monodromy matrix of a symmetric orbit at a symmetric point is a \emph{Wonenburger} matrix:
\begin{equation}\label{symsymp}
M=M_{A,B,C}=\left(\begin{array}{cc}
A & B\\
C & A^T
\end{array}\right)\in Sp(2n),
\end{equation}
where
\begin{equation*}
B=B^T,\quad C=C^T,\quad AB=BA^T,\quad
A^TC=CA,\quad A^2-BC=id,
\end{equation*}
equations which ensure that $M$ is symplectic. The eigenvalues of $M$ are determined by those of the first block $A$\cite{FM}:
\begin{itemize}
    \item If $\lambda$ is an eigenvalue of $M$ then its stability index $a(\lambda)=\frac{1}{2}(\lambda + 1/\lambda)$ is an eigenvalue of $A$. 
    \item If $a$ is an eigenvalue of $A$ then $\lambda(a)=a+\sqrt{a^2-1}$ is an eigenvalue of $M$.
\end{itemize}

\subsection{B-signs} Assume $n=2,3$. Let $x$ be a symmetric orbit with monodromy $M_{A,B,C}$ at a symmetric point. Assume $a$ is a real, simple and nontrivial eigenvalue of $A$ (i.e.,\ $\lambda(a)$ is elliptic or hyperbolic)). Let $v$ be an eigenvector of $A^T$ with eigenvalue $a$, i.e.,\ $A^Tv=a\cdot v$. The \emph{B-sign} of $\lambda(a)$ is 
$$
\epsilon(\lambda(a))=\mbox{sign}(v^TBv)=\pm.
$$

One easily sees that this is independent of $v$, and the basis chosen to write down the monodromy matrix. Note that if $n=2$, we have two $B$-signs $\epsilon_1,\epsilon_2$, one for each symmetric point; and if $n=3$, we have two \emph{pairs} of $B$-signs $(\epsilon^1_1,\epsilon^1_2),(\epsilon^2_1,\epsilon^2_2)$, one for each symmetric point and each eigenvalue. 

The second and fourth authors have recently shown that a planar symmetric orbit is negative hyperbolic iff the $B$-signs of its two symmetric points differ \cite{FM2}. One can define the \emph{$C$-signs} similarly, obtained by replacing the $B$-block, with the $C$-block of $M$, and $A^T$, by $A$.

\subsection{Conley--Zehnder index} 
\begin{figure}
    \centering
    \includegraphics[width=0.8\linewidth]{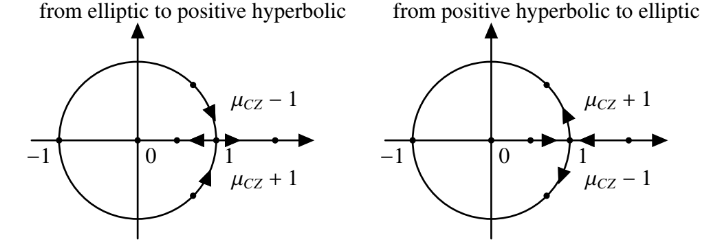}
    \caption{$\mu_{CZ}$ jumps by $\pm 1$ when crossing $1$, according to direction of bifurcation, as shown. If it stays elliptic, the jump is by $\pm 2$. This is determined by the $B$-sign.}
    \label{fig:CZ-jumps}
\end{figure}

The CZ-index is part of the index theory of the symplectic group. It assigns a winding number to non-degenerate orbits. In practical terms, it helps understand which families of orbits connect to which (CZ-index stays constant if no bifurcation occurs, and jumps under bifurcation as shown in Figure~\ref{fig:CZ-jumps}). It may be defined as follows.

\medskip

\textbf{Planar case.} Let $n=2$, $x$ planar orbit with (reduced) monodromy $M^{red}_x$, and $x^k$ its $k$-fold cover which we assume to be non-degenerate for all $k \geq 1$.
\begin{itemize}

\medskip

    \item \textbf{Elliptic case:} $M^{red}_x$ is conjugated to a rotation,
    \begin{equation}\label{reduced_planar}
        M^{red}_x \sim \left(\begin{array}{cc}
       \cos \varphi  & -\sin \varphi \\
    \sin \varphi &  \cos \varphi
    \end{array}\right),
    \end{equation}
with Floquet multipliers $e^{\pm i\varphi}$. Here, $\varphi$ is the rotation angle. Then
    \begin{tcolorbox}\vspace{-0.1cm}
    $$
    \mu_{CZ}(x^k)=1+2\cdot \lfloor k\cdot\varphi/2\pi\rfloor
    $$     
    \end{tcolorbox}
   
 In particular, it is odd, and jumps by $\pm$ 2 if the eigenvalue $1$ is crossed in a family. Recall from (\ref{symsymp}) that for symmetric periodic orbits we have $M^{red}_x = \left(\begin{array}{cc}
       a  & b \\
    c &  a
    \end{array}\right)$. Moreover, in view of (\ref{reduced_planar}) if $b < 0$ then the rotation is determined by $\varphi$ and if $b>0$ then the rotation is determined by $-\varphi$; this determines the CZ-index jump, see Figure \ref{fig:CZ-jumps}.

\medskip
    
    \item \textbf{Hyperbolic case:} $M_x^{red}$ is diagonal up to conjugation, $$
    M^{red}_x \sim \left(\begin{array}{cc}
\lambda  & 0 \\
    0 &  1/\lambda
    \end{array}\right),
    $$
    with Floquet multipliers $\lambda,1/\lambda$. Then
    \begin{tcolorbox}\vspace{-0.1cm}
    $$
    \mu_{CZ}(x^k)=k\cdot n,
    $$    
    \end{tcolorbox}
    
    where $D\phi_t^H$ rotates the eigenspaces by angle $\frac{\pi n t}{T}$, with $n$ even/odd if $x$ positive/negative hyperbolic. Notice that for symmetric periodic orbits the signatures of $b$ and $c$ are equal. 
\end{itemize}

Note that in both cases above, in order to compute the CZ-index via the above formulae, we need to know the linearized flow along the \emph{whole} of the orbit. That is, what matters is the path connecting the identity to the monodromy matrix, obtained by linearizing at any point of the orbit, and performing a full turn around the orbit. In the elliptic case, the rotation angle $\varphi$ is then computed as a real number, and \emph{not} modulo $2\pi$, as it counts the number of rotations of the linearized flow along the whole periodic orbit.

\medskip
    
\textbf{Spatial case.} Let $n=3$. Assume that the reflection along the $xy$-plane gives rise to a symplectic symmetry of $H$ (e.g.,\ the 3BP). If $x\subset \mathbb R^3$ is a planar orbit, then we have a symplectic splitting into planar and spatial blocks
$$
M^{red}_x\sim\left(\begin{array}{cc}
    M_p^{red} &  0\\
    0 & M_s
\end{array}\right)\in Sp(4),\quad M_p^{red}, M_s \in Sp(2).
$$
Then $$\mu_{CZ}(x)=\mu^p_{CZ}(x)+\mu^s_{CZ}(x),$$ where each summand corresponds to $M_p^{red}$ and $M_s$ respectively. We have that
\begin{tcolorbox}
    
\begin{itemize}
    \item Planar to planar bifurcations correspond to jumps in $\mu_{CZ}^p$.
    \item Planar to spatial bifurcations correspond to jumps of $\mu_{CZ}^s$.
\end{itemize}

\end{tcolorbox}

A general definition of the CZ-index will be given in Appendix \ref{app:CZ_index_numerical}, which is needed, e.g.,\ to study spatial-to-spatial bifurcations, and provides a direct way to numerically compute the CZ-indices.  The computations of CZ-indices of families can also be carried out by not directly on the definition, but rather on knowing them analytically for special families (e.g.,\ in the Kepler problem), and then determining the jumps at bifurcations arising after deformation, for which the $B$-signs are necessary, as explained above. This was the approach used by Aydin in \cite{Cengiz}. 

\subsection{Floer numerical invariants}
\begin{figure}
    \centering
    \includegraphics[width=0.42\linewidth]{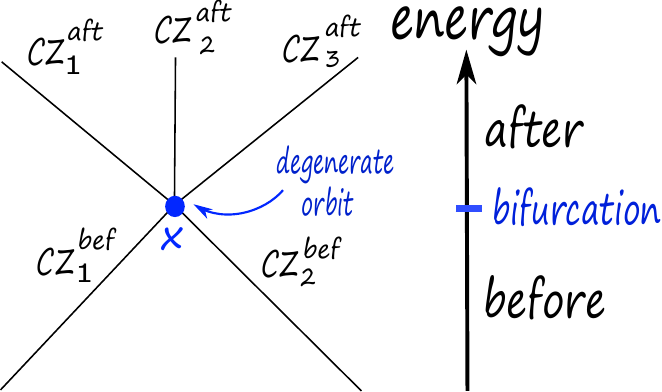}
    \caption{A sketch of a bifurcation at a degenerate orbit, with the before/after orbits determined by the deformation parameter (the energy), each branch with its own CZ-index. The Floer number is a signed count of orbits which stays invariant.}
    \label{fig:bifurcation_sketch}
\end{figure}

Recall that bifurcations occurs when studying families $t\mapsto x_t$ of periodic orbits, as a mechanism by which at some parameter time $t=t_0$ the orbit $x_{t_0}$ becomes degenerate, and several new families may bifurcate out of it; see Figure~\ref{fig:bifurcation_sketch}. The Floer numbers are meant to give a simple test to keep track of all new families. We will first need the following technical definition: a periodic orbit $x$ is \emph{good} if $$\mu_{CZ}(x^k)=\mu_{CZ}(x) \mbox{ mod }2$$ for all $k\geq 1$. Otherwise, it is \emph{bad}. In fact, a planar orbit is bad iff it is an even cover of a negative hyperbolic orbit. And a spatial orbit is bad iff it is an even cover of either an elliptic-negative hyperbolic or a positive-negative hyperbolic orbit. Note that a good planar orbit can be bad if viewed in the spatial problem.

Given a bifurcation at $x$, the \emph{SFT-Euler characteristic} (or the \emph{Floer number}) of $x$ is
\begin{tcolorbox}
$$
\chi(x)=\sum_i (-1)^{CZ^{bef}_i}=\sum_j (-1)^{CZ^{aft}_j}. 
$$    
\end{tcolorbox}
The sum on the LHS is over \textbf{good} orbits \emph{before} bifurcation, and RHS is over \textbf{good} orbits \emph{after} bifurcation. As these numbers only involve the parity of the CZ-index, one has simple formulas which bypass the computation of this index, as they only involve the Floquet multipliers:

\begin{itemize}
    \item \textbf{Planar case.} $\chi(x)=\#\big\{\textrm{ good}\,\,\mathcal H^{+}\big\}-\#\big\{\mathcal E,\,\mathcal{H}^-\big\}.$
\item \textbf{Spatial case.} $\chi(x)=\#\big\{\mathcal H^{--}\,,\mathcal{EH}^-\,,\mathcal E^2\,,\textrm{ good}\,\,\mathcal H^{++}\,,\mathcal N\big\}-\#\big\{\mathcal H^{-+},\,\textrm{good}\,\,\mathcal{EH}^+\big\}.$
\end{itemize}

Here, $\mathcal{E}$ denotes \emph{elliptic}, $\mathcal{H}^\pm$ denotes \emph{positive/negative hyperbolic}, and $\mathcal{N}$ denotes \emph{nonreal} quadruples $\lambda,1/\lambda,\overline{\lambda},1/\overline{\lambda}$. The above simply tells us which type of orbit comes with a plus or a minus sign (the formula should be interpreted as either before or after).

\medskip

\textbf{Invariance.} The fact that the sums agree before and after --\emph{invariance}-- follows from deep results from \emph{Floer theory} in symplectic geometry\footnote{For generic families of Hamiltonians on $4$-dimensional phase spaces this can alternatively be proved by the using the normal forms of Meyer \cite{Meyer1970}. See for instance the appendix in \cite{FKM}.}. We will accept this as a fact, and use it as follows:

\medskip

\begin{tcolorbox}
    The Floer number can be used as a \textbf{test}: if the sums do \emph{not} agree, we know the algorithm missed an orbit.
\end{tcolorbox}

The invariant above works for arbitrary periodic orbits. There is a similar Floer invariant for \emph{symmetric} orbits \cite{FKM}.

\subsection{Global topological methods}

These methods encode: bifurcations; stability; eigenvalue configurations; obstructions to existence of regular families; and $B$-signs, in a visual and resource-efficient way. The main tool is the \emph{GIT sequence} \cite{FM}, a refinement of the Broucke stability diagram via implementing the $B$-signs. This is a sequence of three branched spaces (or \emph{layers}), together with two maps between them, which collapse certain branches together. Each branch is labeled by the $B$-signs. A symmetric orbit gives a point in the top layer, and an arbitrary orbit, in the middle layer. The base layer is $\mathbb R^n$ (the space of coefficients of the characteristic polynomial of the first block of $M_{A,B,C}$). Then a family of orbits gives a path in these spaces, so that their topology encodes valuable information. The details are as follows. 

\medskip

\textbf{GIT sequence: 2D.} Let $n=2$, $\lambda$ eigenvalue of $M^{red}\in Sp(2)$, with stability index $a(\lambda)=\frac{1}{2}(\lambda+1/\lambda)$. Then $\lambda=\pm 1$ iff $a(\lambda)=\pm 1$; $\lambda$ positive hyperbolic iff $a(\lambda)>1$; $\lambda$ negative hyperbolic iff $a(\lambda)<-1$; and $\lambda$ elliptic (stable) iff $-1<a(\lambda)<1$. The Broucke stability diagram is then simply the real line, split into three components; see Figure~\ref{fig:GIT_sequence}. If two orbits lie in different components of the diagram, then one should expect bifurcations in any family joining them, as the topology of the diagram implies that any path between them has to cross the $\pm 1$ eigenvalues.

One can think that the stability index ``collapses'' the two elliptic branches in the middle layer of Figure \ref{fig:GIT_sequence} together. These two branches are distinguished by the $B$-signs, coinciding with the Krein signs \cite{Kre2,Kre3}. There is an extra top layer for symmetric orbits, where now each hyperbolic branch separates into two, and there is a collapsing map from the top to middle layer. Note that to go from one branch to the other, the topology of the layer implies that the eigenvalue 1 needs to be crossed. This means that one should expect bifurcations in any (symmetric) family joining them, \emph{even if} they project to the same component of the Broucke diagram. To sum up:

\begin{figure}[ht]
    \centering
    \includegraphics[width=0.7\linewidth]{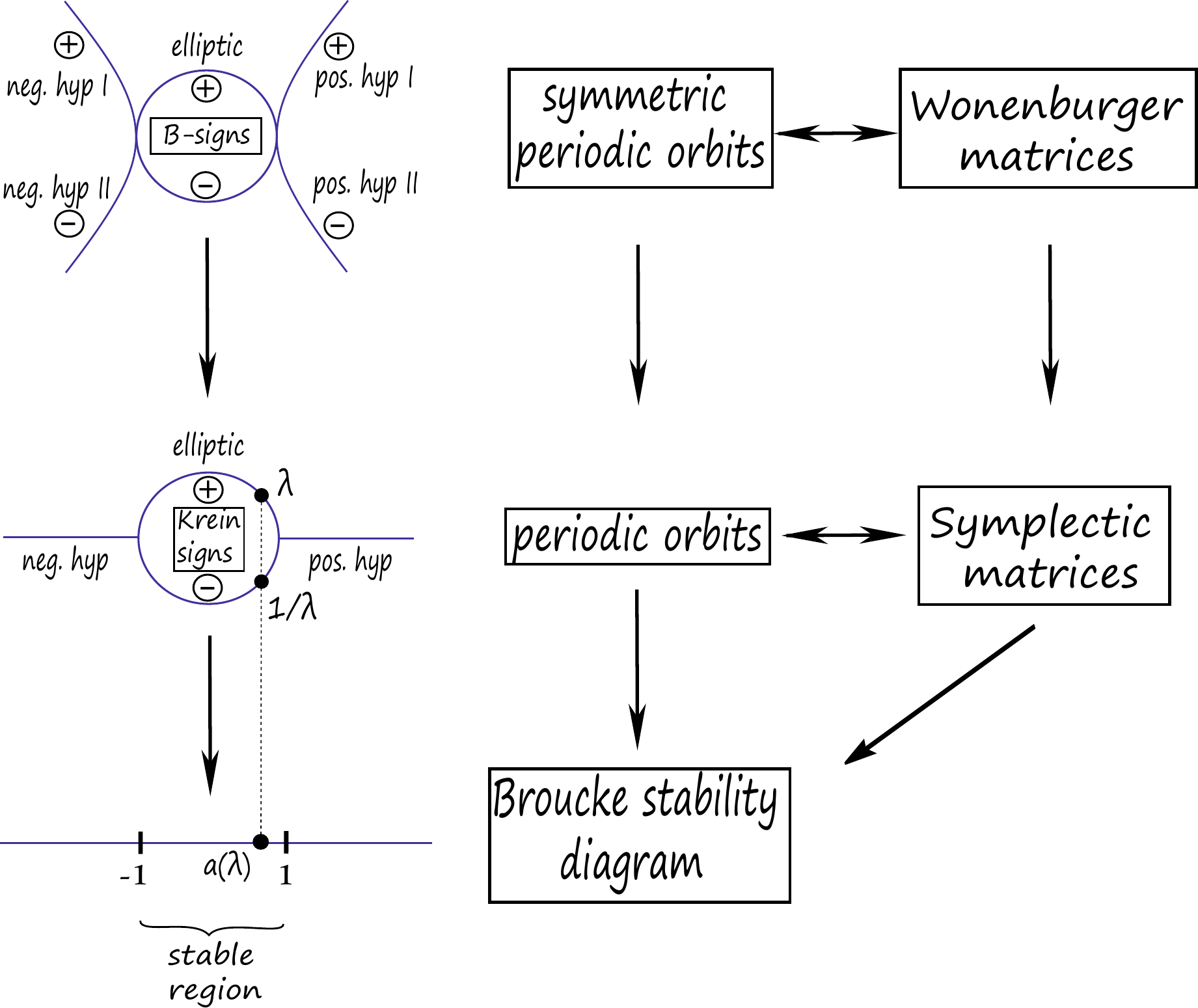}
    \caption{The 2D GIT sequence. One obtains more refined information for symmetric orbits.}
    \label{fig:GIT_sequence}
\end{figure}

\medskip

\begin{tcolorbox}
\begin{itemize}
    \item $B$-signs ``separate'' hyperbolic branches, for symmetric orbits.

\medskip
    
    \item If two points lie in different components of the Broucke diagram, one should expect bifurcation in any path joining them.

\medskip

    \item If two points lie in the same component of the Broucke diagram, but if $B$-signs differ, one should \emph{also} expect bifurcation in any path joining them. 
    \end{itemize}
\end{tcolorbox}

\begin{figure}
    \centering
    \includegraphics[width=0.55\linewidth]{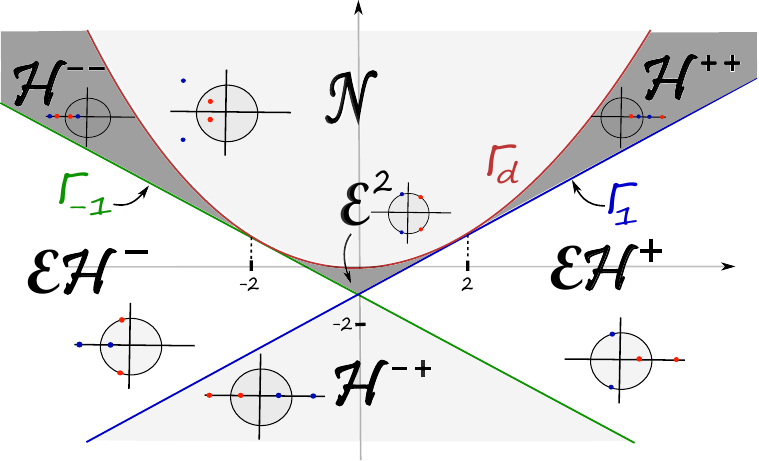}
    \caption{The 3D Broucke stability diagram. Here, $\Gamma_{\pm 1}$ corresponds to eigenvalue $\pm 1$, $\Gamma_d$ to double eigenvalue, $\mathcal{E}^2$ to doubly elliptic (\textbf{stable region}), and so on \cite{FM}.}
    \label{fig:Broucke_3D}
\end{figure}

\begin{figure}
    \centering
    \includegraphics[width=1\linewidth]{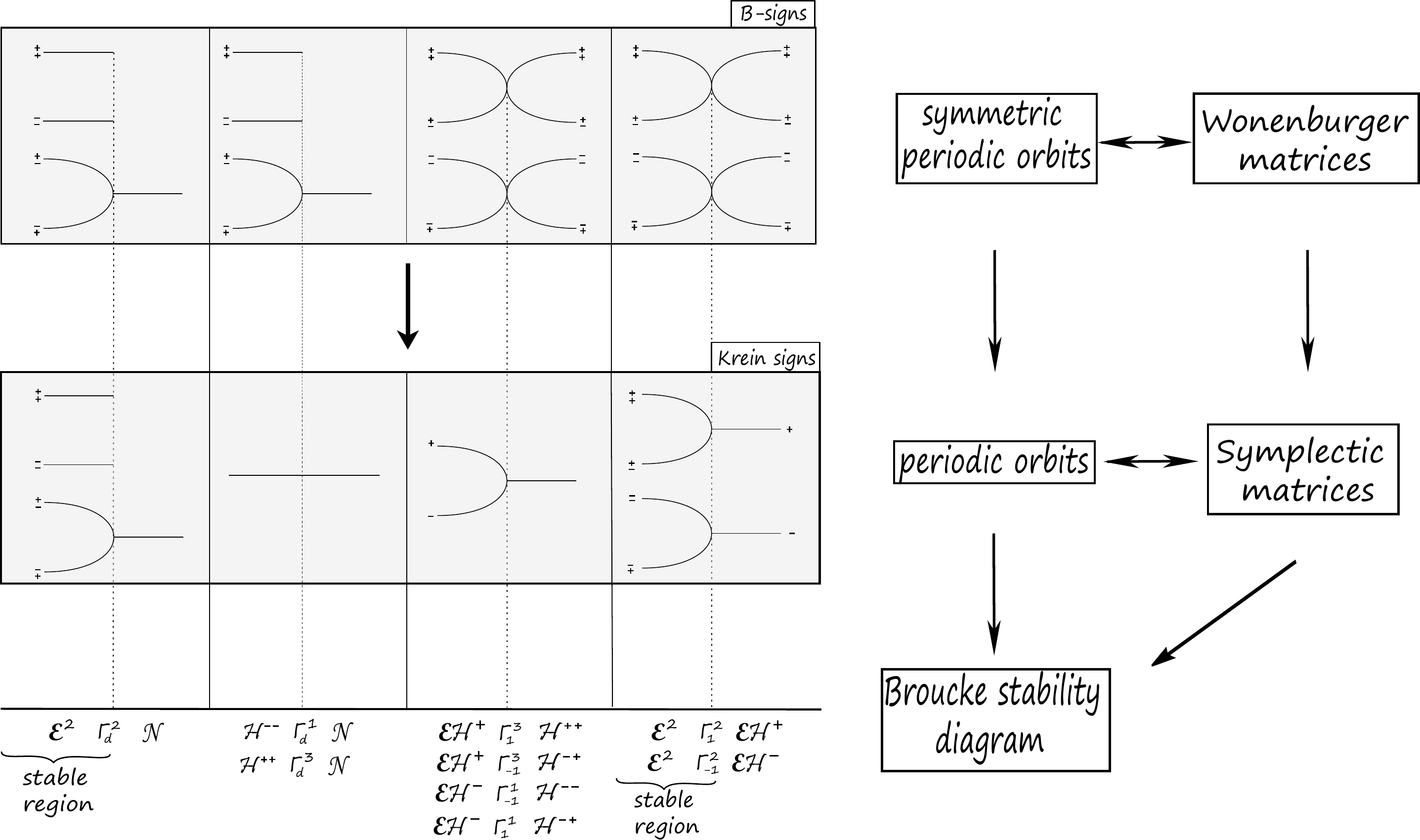}
    \caption{The branches (represented as lines) are two-dimensional, and come together at the 1-dimensional ``branching locus'' (represented as points), where we cross from one region to another of the Broucke diagram.}
    \label{fig:GIT_sequence_3D}
\end{figure}

\textbf{GIT sequence: 3D.} Let $n=3$. Given $M^{red}=M_{A,B,C}\in Sp(4)$, its \emph{stability point} is $p=(\mbox{tr}(A),\det(A)) \in \mathbb R^2$. The plane splits into regions corresponding to the eigenvalue configuration of $M^{red}$, as in Figure~\ref{fig:Broucke_3D}. The GIT sequence \cite{FM} adds two layers to this diagram, as shown in Figure~\ref{fig:GIT_sequence_3D}. The top layer has two extra branches than the middle one, for each hyperbolic eigenvalue.

\medskip

\textbf{Bifurcations in the Broucke diagram.} An orbit family $t\mapsto x_t$ gives a path $t\mapsto p_t\in \mathbb R^2$ of stability points. The family bifurcates if $p_t$ crosses $\Gamma_1$. More generally, let $\Gamma^e_{\varphi}$ be the line with slope $\cos(2\pi \varphi)\in [-1,1]$ tangent to $\Gamma_d=\{y=x^2/4\}$, corresponding to matrices with eigenvalue $e^{2\pi i\varphi}$; and $\Gamma^h_\lambda$ the tangent line with slope $a(\lambda)\in \mathbb R\backslash [-1,1]$, corresponding to matrices with eigenvalue $\lambda$.
    
\begin{figure}
    \centering
\includegraphics[width=0.55\linewidth]{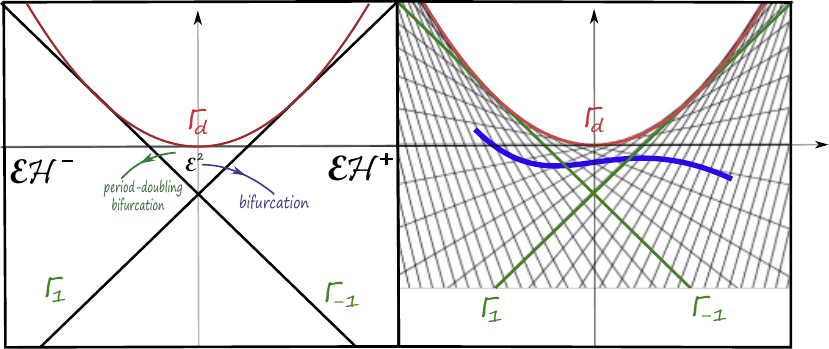}
\caption{Bifurcations are encoded by a pencil of lines.}
   \label{fig:complete_bifurcation}
\end{figure}

\begin{tcolorbox}
A $k$-fold bifurcation happens when crossing $\Gamma^e_{l/k}$ for some $l$.    
\end{tcolorbox}

That is, higher order bifurcations are encoded by a pencil of lines tangent to a parabola, as in Figure~\ref{fig:complete_bifurcation}.

\medskip

\begin{figure}
        \centering
        \includegraphics[width=1\linewidth]{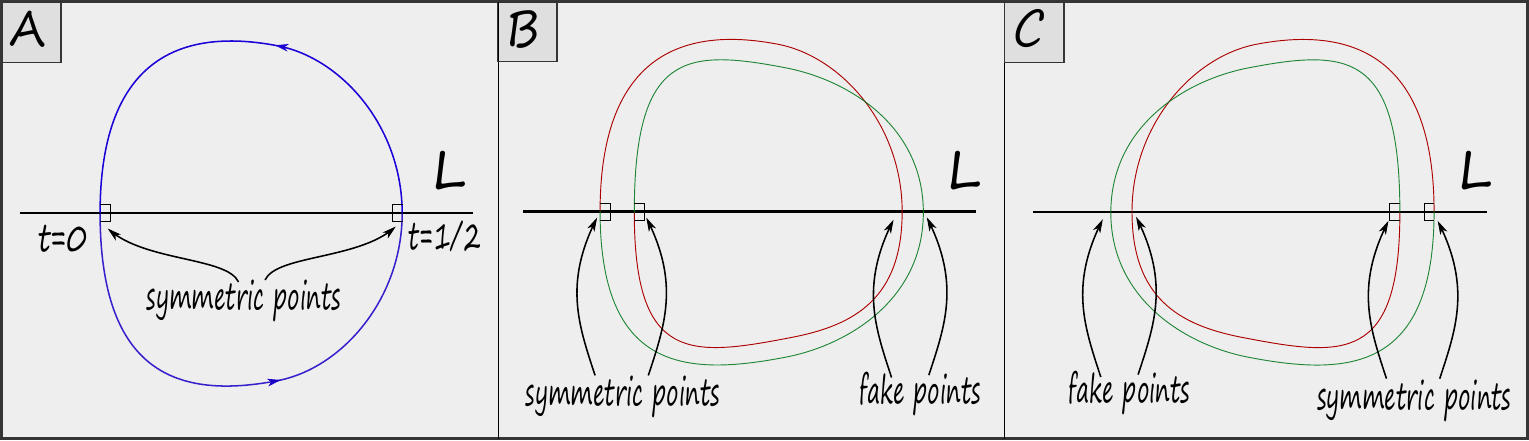}
        \caption{Symmetric period doubling bifurcation. The \emph{fake} symmetric points, while close to intersection points, do \emph{not} intersect the fixed-point loci.}
        \label{fig:period_doubling}
    \end{figure}

\textbf{Example: symmetric period doubling bifurcation.} We finish this section with an example where our invariants give new information. Consider a symmetric orbit $x$ going from elliptic to negative hyperbolic. A priori there could be two bifurcations, one for each symmetric point (B or C in Figure~\ref{fig:period_doubling}). However, invariance of $\chi(x^2)$ implies only \emph{one} can happen (note $x^2$ is \emph{bad}). And where the bifurcation happens is determined by the $B$-sign, occurring at the symmetric point in which the $B$-sign does \emph{not} jump; or alternatively, where the $C$-sign jumps.

\subsection{Circular restricted three-body problem}

\begin{figure}[t]
	\centering\includegraphics[width=2.2in]{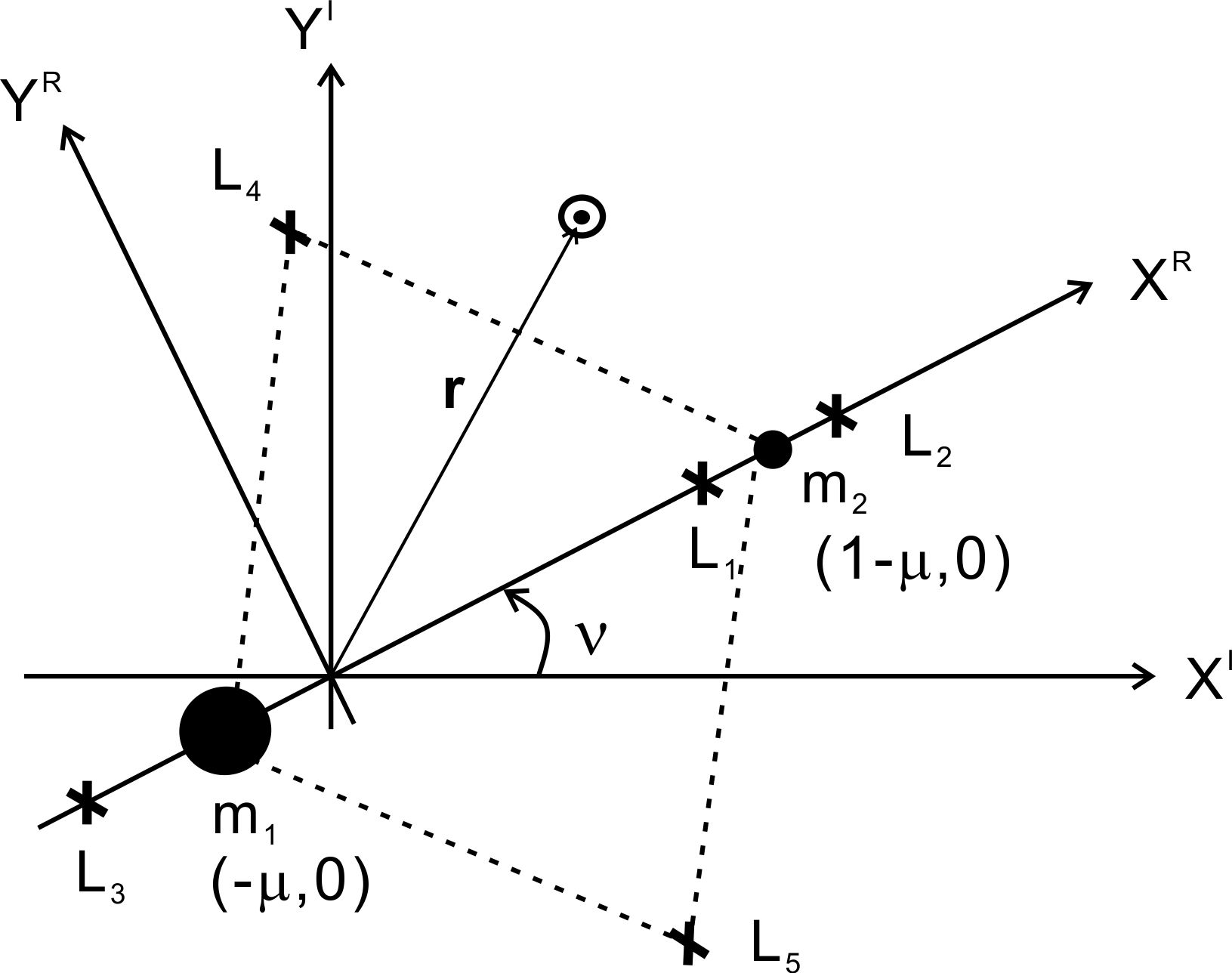}
	\caption{A schematic CRTBP configuration showing $x_1=m_1$, $x_2=m_2$, and two of the libration points in a non-dimensional rotating coordinate system $X^R-Y^R$,~ $Z^R$($Z^I$) are in the out-of-plane direction.}
	\label{fig:crtbp}
\end{figure}

The Circular Restricted Three-Body Problem (CRTBP) shown in Figure~\ref{fig:crtbp} describes the motion of an infinitesimal mass with two primaries under mutual gravitational attraction. A dimensionless rotating coordinate system ($X^R-Y^R-Z^R$) is defined at the barycenter of the two primaries with respect to the inertial frame ($X^I-Y^I-Z^I$), rotating about $Z^I$ with true anomaly $\nu$. 
 % Fig crtbp
The $X$-axis of the rotating coordinate system is aligned with the vector from the larger primary body ($m_1$) to the second primary body ($m_2$). The $Z$-axis is perpendicular to the primaries' orbital plane, and the $Y$-axis completes the right-handed coordinate system. The position vector $\mathbf{r}$ points from the barycenter to the spacecraft in the rotating frame. The non-dimensional mass of the second primary is defined as 
\begin{equation*}
	\mu=\frac{m_2}{m_1+m_2}=m_2,
\end{equation*}
and then the larger body's mass is
\begin{equation*}
	1-\mu=\frac{m_1}{m_1+m_2}=m_1.
\end{equation*}
Define the unit of time so that the mean motion of the primary orbit is $1$. Then the equations of motion for the infinitesimal mass is written as 
\begin{equation*}
\begin{aligned}
	\label{eq: CRTBP}
	\ddot{x}	&=2\dot{y}+x-(1-\mu)\frac{x+\mu}{r_{1}^{3}}-\mu\frac{x-1+\mu}{r_{2}^{3}}\\
	\ddot{y}	&=-2\dot{x}+y-(1-\mu)\frac{y}{r_{1}^{3}}-\mu\frac{y}{r_{2}^{3}}\\
	\ddot{z}	&=-(1-\mu)\frac{z}{r_{1}^{3}}-\mu\frac{z}{r_{2}^{3}}
\end{aligned} 
\end{equation*}
where $r_1^2 = (x+\mu)^2+y^2+z^2$, $r_2^2=(x-1+\mu)^2+y^2+z^2$. No closed form general solution is possible for the model. 

The Hamiltonian describing the CRTBP is given by 
$$
H: (\mathbb R^3\setminus \{M,P\})\times \mathbb R^3\rightarrow \mathbb R,
$$
$$
H(q,p)=\frac{1}{2}\Vert p\Vert^2 - \frac{\mu}{\Vert q- M\Vert } - \frac{1-\mu}{\Vert q- P\Vert } +p_1q_2-p_2q_1, $$
where $q=(q_1,q_2,q_3)$ is the position of a satellite, $p=(p_1,p_2,p_3)$ is its momentum, the mass of the secondary body $m_2$ is fixed at $M=(1-\mu,0,0)$, and the mass of the primary body $m_1$ is fixed at $P=(-\mu,0,0)$. The Jacobi constant $\Gamma$ is then defined by the convention $\Gamma:=-2H$. The Hamiltonian $H$ is invariant under the anti-symplectic involutions 
$$
\rho: (q_1,q_2,q_3,p_1,p_2,p_3) \mapsto (q_1,-q_2,-q_3,-p_1,p_2,p_3),
$$
$$
\widetilde\rho: (q_1,q_2,q_3,p_1,p_2,p_3) \mapsto (q_1,-q_2,q_3,-p_1,p_2,-p_3),
$$
with corresponding fixed-point loci given by 
$$
L=\mbox{Fix}(\rho)=\{q_2=q_3=p_1=0\},
$$ 
$$
\widetilde L=\mbox{Fix}(\widetilde \rho)=\{q_2=p_1=p_3=0\}.
$$
These correspond respectively to $\pi$-rotation around the $x$-axis, and reflection along the $xz$-plane. Their composition $\sigma=\rho \circ \tilde \rho$ is a \emph{symplectic} symmetry corresponding to reflection along the $xy$-plane.

For instance, the \emph{Jupiter-Europa system} then corresponds to a CRTBP with mass ratio $\mu=2.5266448850435e^{-05}$, and the \emph{Saturn-Enceladus system}, to $\mu=1.9002485658670e^{-07}$. 
This information as well as other orbital data can be found in \cite{Jaumann2009} as well as on the webpage
\begin{verbatim}
    https://ssd.jpl.nasa.gov/sats/
\end{verbatim}
We shall use this information below.

\subsection{Hill's lunar problem}

Hill's lunar problem is a limit case of the restricted three-body problem where the infinitesimal mass is assumed very close to the small primary. This problem can therefore be viewed as an approximation to
the Saturn--Enceladus and Jupiter--Europa system, when one lets the mass of Europa go to zero. The Hamiltonian describing the system is
$$
E: (\mathbb R^3\backslash\{0\})\times \mathbb R^3\rightarrow \mathbb R,
$$
$$
E(q,p)=\frac{1}{2}\Vert p \Vert^2-\frac{1}{\Vert q \Vert}+p_1q_2-p_2q_1-q_1^2+\frac{1}{2}q_2^2+\frac{1}{2}q_3^2.
$$
The linear symmetries of this problem have been completely characterized \cite{Cengiz2}. While the planar restricted three-body problem is invariant
under reflection at the $x$-axis, the planar Hill lunar problem is additionally invariant under reflection at the $y$-axis. For the spatial lunar problem, there are more symmetries: $\rho,\widetilde{\rho}$ (which extend the reflection at the $x$-axis), and two additional symmetries $\kappa,\widetilde{\kappa}$ ($\pi$-rotation along the $y$-axis, and reflection along the $yz$-plane; both extend the reflection along the $y$-axis). Their composition is also $\sigma=\kappa \circ \widetilde \kappa$. 

\section{Numerical work}

\subsection{Result I. Planar direct/prograde orbits}

\begin{figure}[hbtp]
    \centering
    \includegraphics{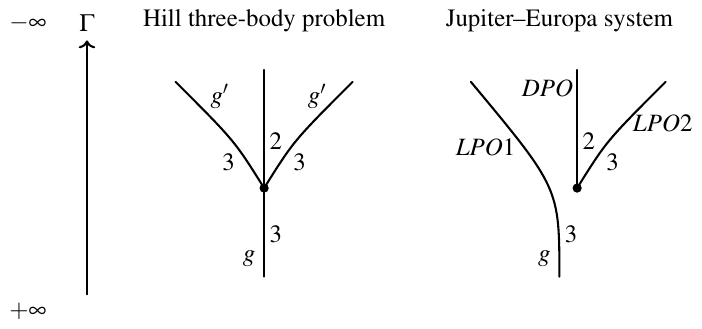}
    \caption{Bifurcation graphs for the planar direct/prograde orbits with CZ-index.}
    \label{fig:pitchfork}
\end{figure}

\begin{figure}[h]
    \centering
    \includegraphics[width=0.7\linewidth]{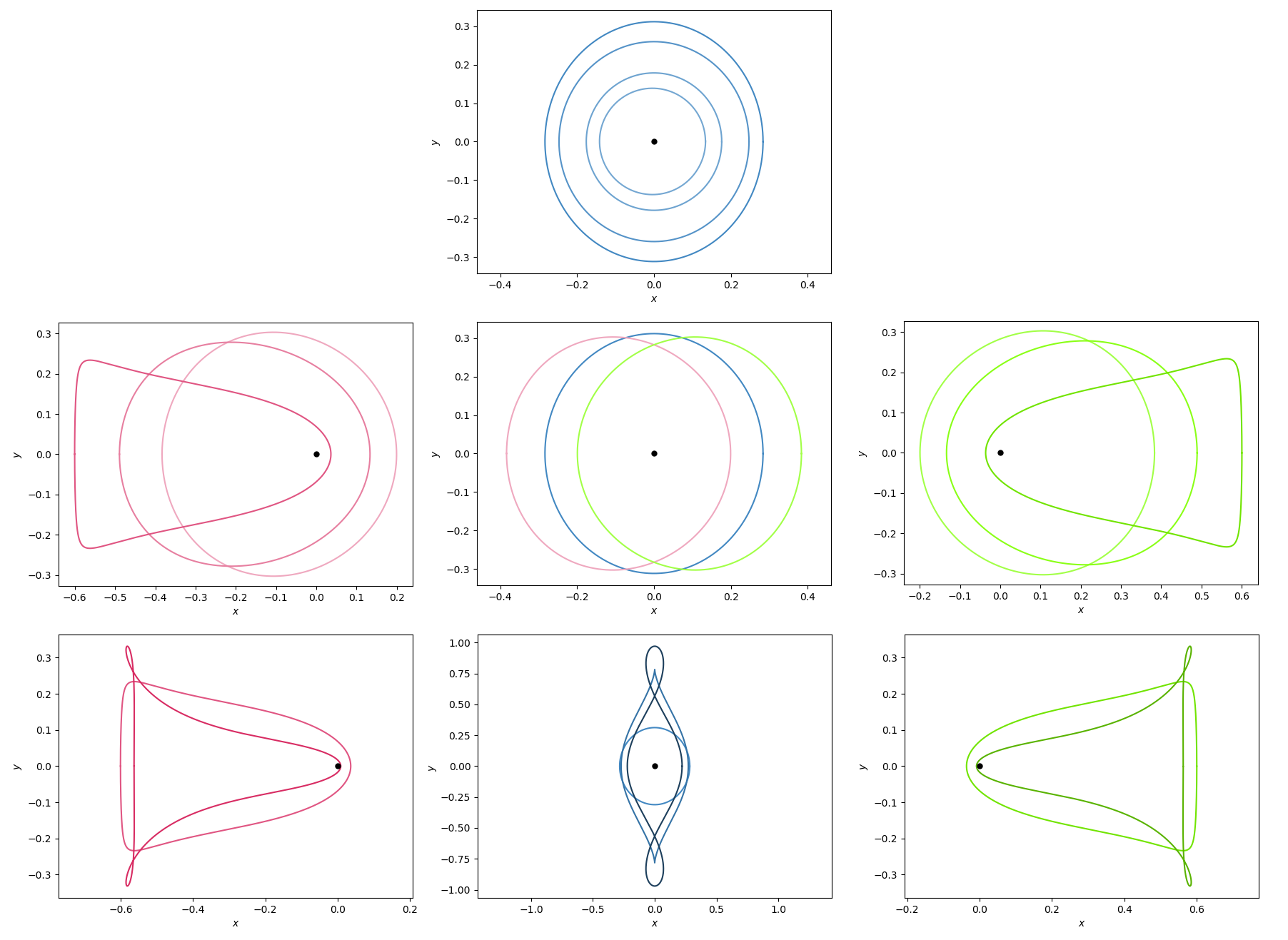}
    \caption{Middle: the $g$ branch whose orbits are in blue and doubly-symmetric w.r.t.\ the $x$- and $y$-axis. From the orbit in the center there bifurcate the two $g'$ branches whose orbits are simply-symmetric w.r.t. the $x$-axis; one $g'$ branch is on the right (the green orbits) and one $g'$ branch is on the left (the purple orbits). The energy increases from top to bottom and from light to dark in each plot.}
    \label{fig:g_g'_hill}
\end{figure}

H\'enon \cite{henon} describes a family $g$ of planar direct periodic orbits which are invariant with respect to both reflections at the $x$ and $y$-axis. This family undergoes a non-generic pitchfork bifurcation, going from elliptic to positive hyperbolic, and where two new families of elliptic orbits, called $g'$, appear; see the plots in Figure \ref{fig:g_g'_hill}. These new families are still invariant under reflection at the $x$-axis, but not under reflection at the $y$-axis. Reflection at the $y$-axis maps one branch of the $g'$-family to the other branch. Figure \ref{fig:pitchfork} shows the bifurcation graph which is constructed as follows: each vertex denotes a degenerate orbit at which bifurcation happens and each edge represents a family of orbits with varying energy, labeled by the corresponding CZ-index. From this data, it is easy to determine the associated Floer number. For instance in Figure \ref{fig:pitchfork} on the left, the Floer number is $(-1)^3 = -1$ before bifurcation, and $(-1)^2 + 2(-1)^3 = -1$ after bifurcation; they coincide, as they should.

\begin{figure}[hbtp]
    \centering
    \includegraphics[width=0.7\linewidth]{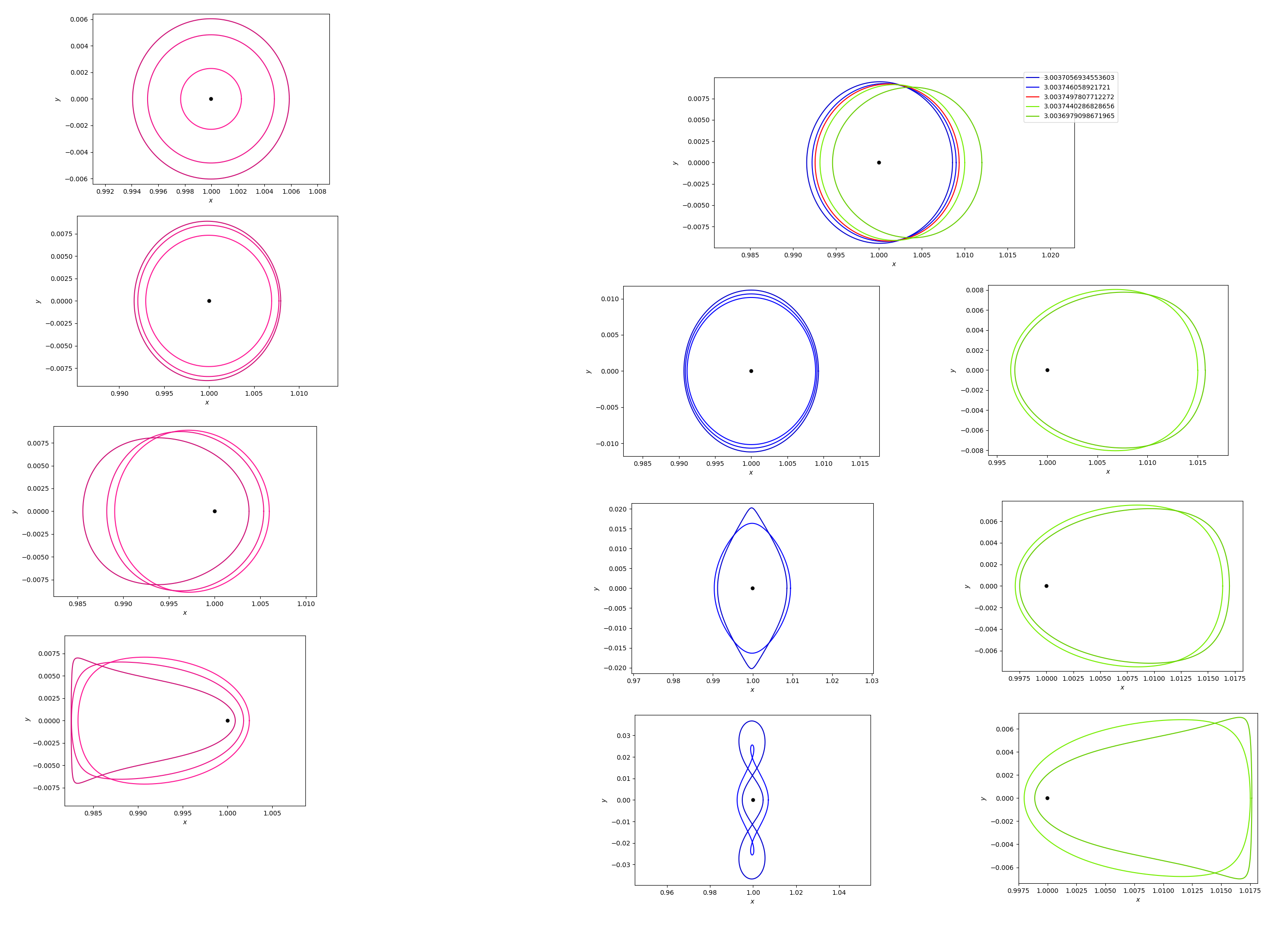}
    \caption{Left: the g-LPO1 branch, where the energy increases from top to bottom. Right: the DPO-LPO2 branch, split into the DPO sub-branch (left, where orbits are planar positive hyperbolic) and of the LPO2 sub-branch (right, where orbits are planar elliptic). The red orbit is the degenerate orbit, undergoing birth-death bifurcation. The Jacobi constant $\Gamma=-2c$ is shown on the upper right side, reaching a maxima at the red orbit.}
    \label{fig:branches}
\end{figure}

By deforming the mass parameter, we may go from Hill's lunar problem to the Jupiter--Europa system; see Figure \ref{fig:pitchfork}. The pitchfork bifurcation deforms to a \emph{generic} situation, where one of the $g'$ branches glues to the before-bifurcation part of the $g$ branch, the result of which we call the \emph{g-LPO1} branch, and where the other $g'$ branch glues to the after-bifurcation part of the $g$ branch, which we call the \emph{DPO-LPO2} branch (undergoing birth-death bifurcation). The $DPO$-orbits are planar positive hyperbolic and the $LPO2$-orbits are planar elliptic. As the symmetry with respect to the $y$-axis is lost, the new orbits will be \emph{approximately} symmetric with respect to the $y$-axis, but not exactly symmetric; similarly, the $y$-symmetric relation between the $g'$ branches persists only approximately for the corresponding deformed orbits. These families are plotted in Figure \ref{fig:branches}, where this behavior is manifest. The data for each new branch is given in Tables \ref{data_g_lpo1}, \ref{data_dpo} and \ref{data_lpo2} in Appendix \ref{app:tables}. Via this bifurcation analysis, one may predict the existence of the DPO-LPO2 branch, which a priori is not straightforward to find. While these families are already known and appear e.g.,\ in page 12 of \cite{database}, this suggests a general mechanism which we will exploit, cf.\ Figure \ref{fig:bifurc_diagram}, and Figure \ref{fig:bifurc_diagram_new}. Note that \cite{database} provides an online data base for \emph{planar} and \emph{$x$-axis symmetric} periodic orbits, and we match their notation for orbits (DPO, LPO, etc.). The novelty of this article is to focus on \emph{spatial} bifurcations of these planar orbits, employing our novel methods and tools.

\begin{figure}[h]
    \centering
    \includegraphics[width=0.9\linewidth]{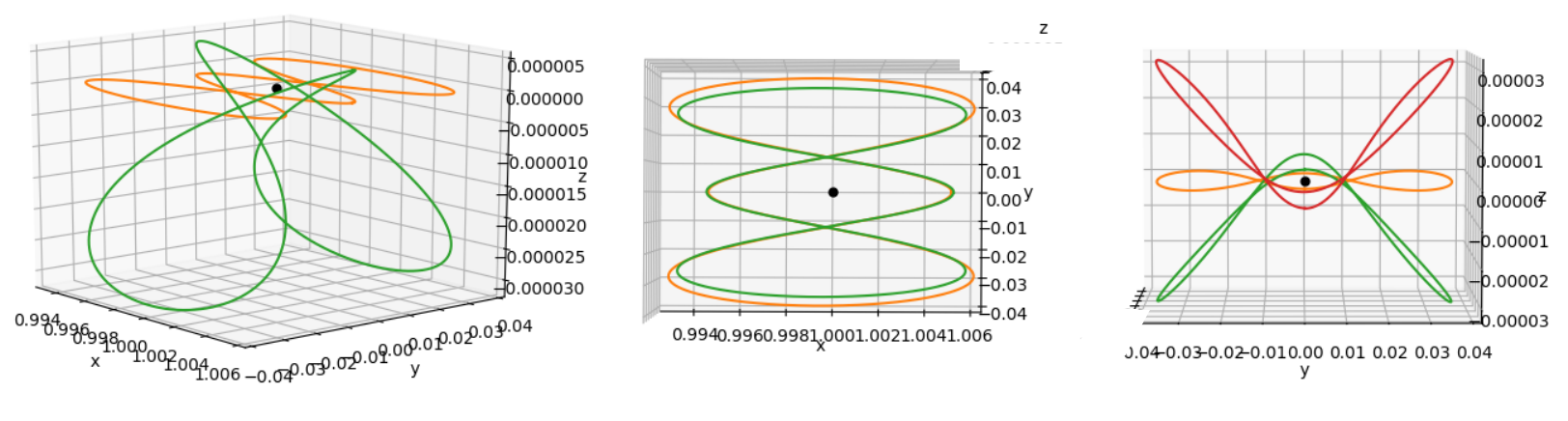}
    \caption{Jupiter-Europa: A planar-to-spatial bifurcation of a simple closed planar DPO orbit, from the side, from above and with its symmetric family by using the reflection at the $xy$-plane.}
    \label{fig:DPO}
\end{figure}

\subsection{Result II. Bifurcation graphs with the same topology}

\begin{figure}[h]
	\centering
	\includegraphics[]{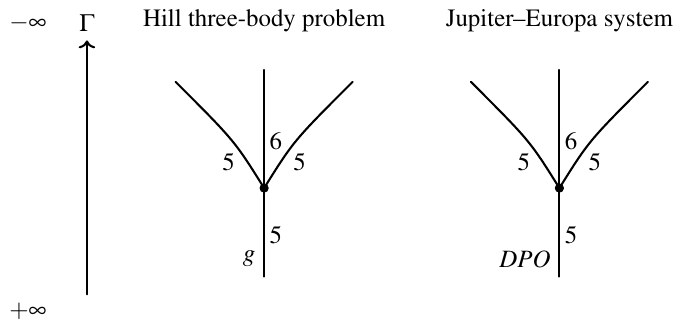}
	\caption{Left:\ The bifurcation graph between simple closed $g$-orbit and the new families of spatial orbits generated by the spatial index jump in Hill's system. Right: In the Jupiter--Europa system. The horizontal symmetry corresponds to the reflection at the $xy$-plane.}
	\label{bifurcation_graph_g_spatial_index_jump}
\end{figure}

In the Jupiter--Europa system, the spatial CZ-index of the simple closed DPO-orbit at around $\Gamma = 3.00109352$ jumps by $+1$, see Table \ref{data_dpo} in Appendix \ref{app:tables}. Therefore it generates a planar-to-spatial bifurcation, see the plot in Figure \ref{fig:DPO}. As in Hill's problem, this new family of spatial orbits appears twice by using the reflection at the $xy$-plane. Surprisingly, compared to Figure \ref{fig:pitchfork}, because the symmetry is preserved, the bifurcation graph has the same topology after deformation and is still non-generic, see the graph in Figure \ref{bifurcation_graph_g_spatial_index_jump}.

\subsection{Result III. Bifurcation graphs between prograde and retrograde orbits}

\begin{figure}[h]
    \centering
    \includegraphics[width=0.8\linewidth]{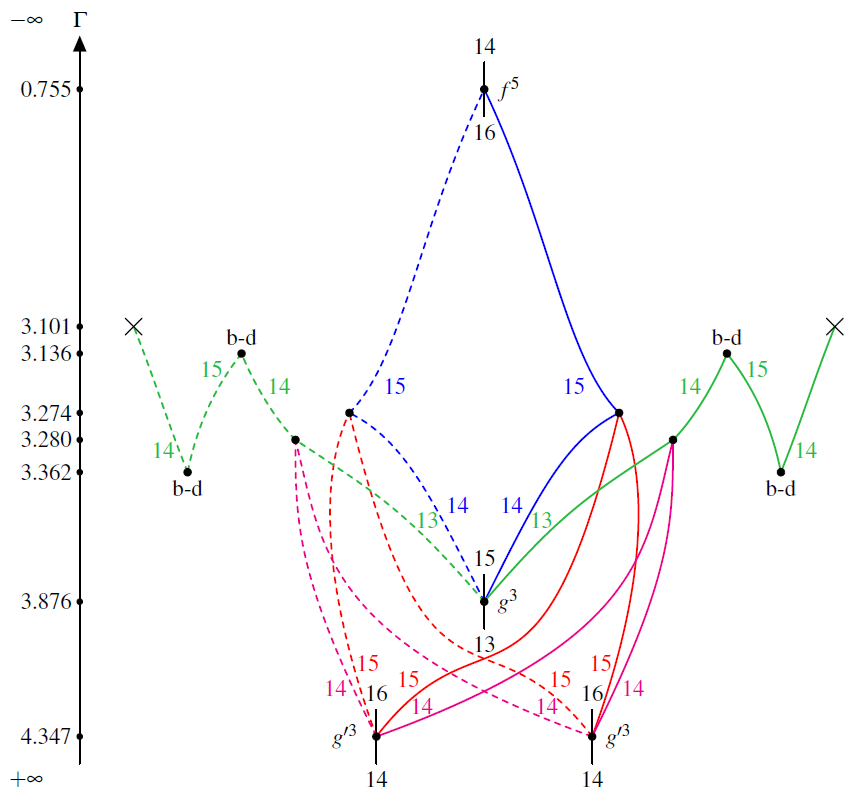}
    \caption{Bifurcation graph for Hill's lunar problem by the first author \cite{Cengiz}, between the 3rd cover of $g$, the 3rd cover of $g'$ and the 5th cover of $f$, based on work of Kalantonis \cite{Kalantonis}. A cross means collision, and b-d means birth-death bifurcation. The horizontal symmetry in the diagram, relating full and dashed edges, means that the corresponding families are related by a symmetry. For instance, the non-dashed red 15 on the right is related by the dashed red 15 on the right by reflection along the xy-plane. The other red 15 families on the left are obtained by applying the extra two spatial symmetries $\kappa,\widetilde \kappa$. Similarly for the pink 14 families. The blue and green families are doubly symmetric; one of the symmetries breaks at bifurcation, where the red and pink families appear.}
    \label{fig:bifurc_diagram}
\end{figure}

\begin{figure}[t]
    \centering
    \includegraphics[width=0.48\linewidth]{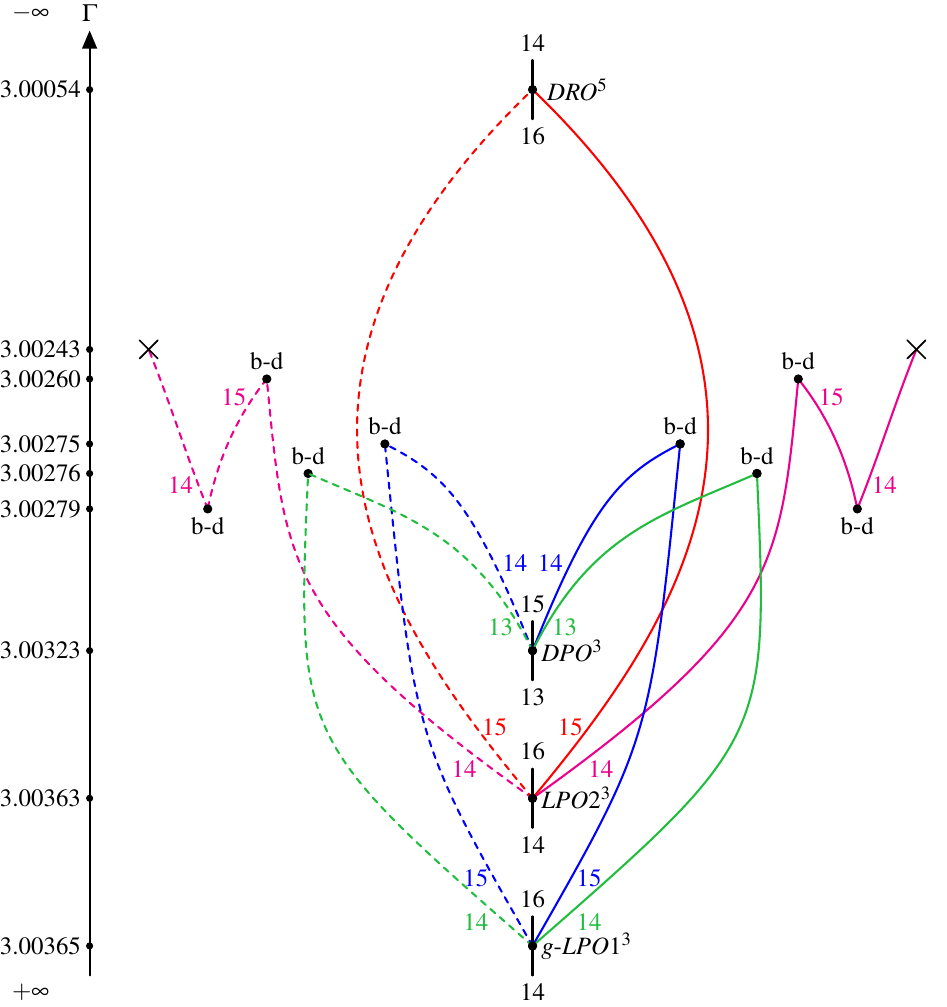}
    \caption{Bifurcation graph for the Jupiter--Europa system, between $g$-$LPO1^3$, $DPO^3$,  $LPO2^3$, and $DRO^5$. The data for the pink family is collected in Table \ref{data_lpo2_3rd_cover_1}, for the red one in Table \ref{data_lpo2_3rd_cover_2}, for the blue one in Table \ref{data_lpo1_3rd_cover_1} and for the green in Table \ref{data_lpo1_3rd_cover_2}. Some orbits of the blue and green family are plotted in Figure \ref{fig:LPO1_3rd}. The horizontal symmetry is reflection along the $xy$-plane. Note that non-dashed red 15 and non-dashed blue 15 are no longer related by a symmetry.}
    \label{fig:bifurc_diagram_new}
\end{figure}

A bifurcation graph relating third covers of $g, g'$, and fifth covers of planar retrograde orbits, known as family $f$, was obtained by the first author \cite{Cengiz}; see Figure \ref{fig:bifurc_diagram}. The third covers of LPO2 and fifth covers of DRO were found using Cell-Mapping \cite{Koh2021}. Taking Figure \ref{fig:bifurc_diagram} as a starting point, we compare it to the Jupiter--Europa system. The result is plotted in Figure \ref{fig:bifurc_diagram_new}. 

Let us focus on the two unlabeled vertices on the right of Figure \ref{fig:bifurc_diagram} which are not of birth-death type. After deformation, the (red) family starting at $g'^3$ on the right of CZ-index 15 glues to the (blue) family of the same index ending in $f^5$, resolving the vertex at which they meet; note that similarly as in Result I, $f$ is replaced with DRO, and $g'$, with LPO2. The two other families meeting at the same vertex coming from $g'^3$ and $g^3$ now glue to a family undergoing birth-death, where now $g'$ is replaced by $g$-$LPO1$, and $g$, with $DPO$. A similar phenomenon happens at the other vertex, where the (pink) family starting at $g'^3$ with CZ-index 14 on the right glues to the (green) family of the same index, and the other two families now undergo birth-death. These families might have been hard to find without this analysis. %As far as we know, these two families (blue and green in Figure \ref{fig:bifurc_diagram_new}) appear to be new.
 
\begin{figure}
    \centering
    \includegraphics[width=0.8\linewidth]{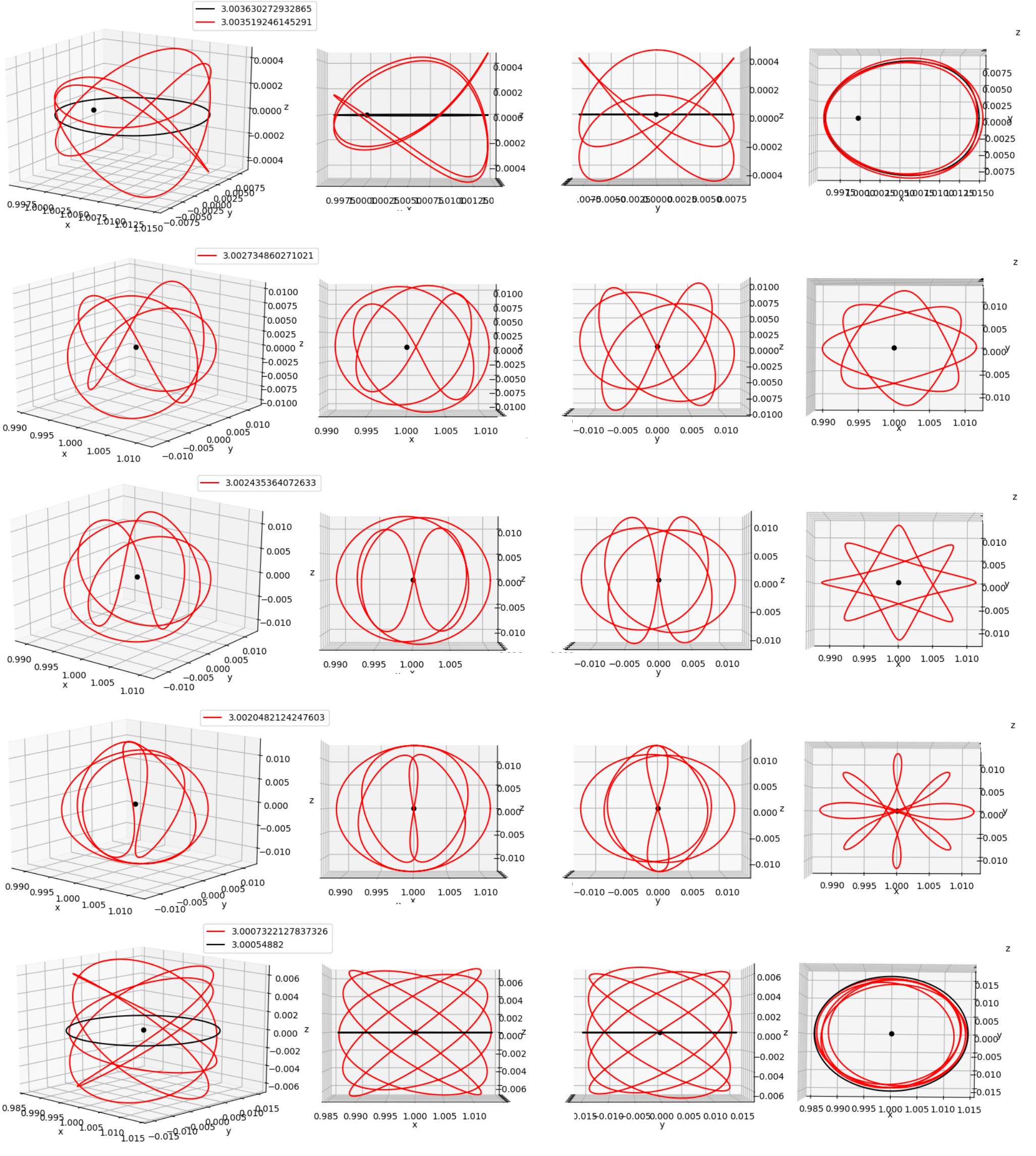}
    \caption{The red prograde to retrograde spatial connection, with CZ-index 15, for varying energy. The black planar orbit is an LPO2, and the red family bifurcates from its third cover. A row corresponds to the same orbit from different angles.}
    \label{fig:pro-ret}
\end{figure}

Another notable feature is the (red) family between $LPO2^3$ and $DRO^5$ of CZ-index 15. This is a \emph{spatial} family connecting two \emph{planar} orbits, one of which is retrograde ($DRO^5$), and the other, prograde ($LPO2^3$). This family is plotted in Figure \ref{fig:pro-ret}.

\begin{figure}[b]
    \centering
    \includegraphics[width=0.65\linewidth]{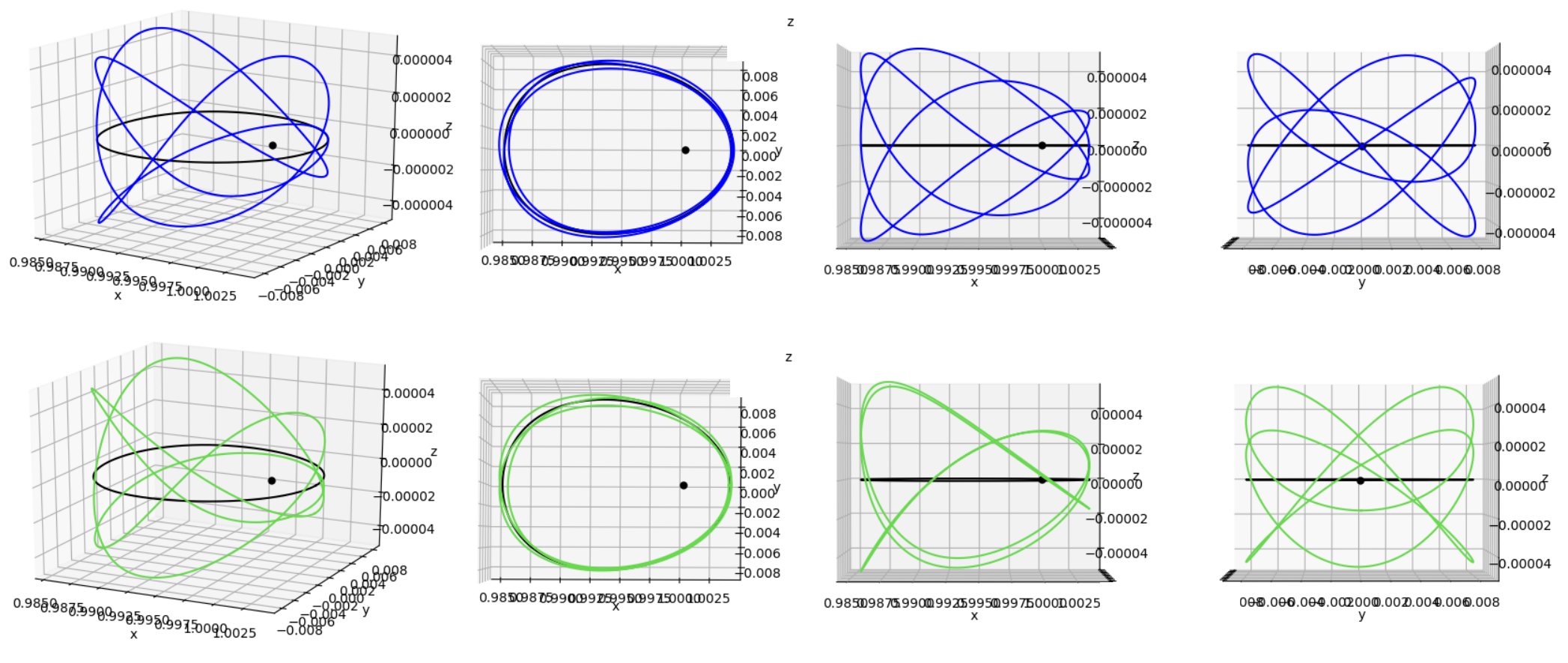}
    \caption{Jupiter-Europa: Two families of spatial orbits branching out from the $g$-$LPO1^3$ orbit; above: these orbits are symmetric w.r.t. the $x$-axis and their data is collected in Table \ref{data_lpo1_3rd_cover_1}. This is the (blue) family of CZ-index 15 in Figure \ref{fig:bifurc_diagram_new}; below: these orbits are symmetric w.r.t.\ the $xz$-plane and their data is collected in Table \ref{data_lpo1_3rd_cover_2}. This is the (green) family of CZ-index 14 in Figure \ref{fig:bifurc_diagram_new}. %Each family has a symmetric family by using the reflection at the ecliptic. 
    }
    \label{fig:LPO1_3rd}
\end{figure}

\subsection{Result IV. Bifurcation graph for Saturn--Enceladus}
\begin{figure}
    \centering
    \includegraphics[width=0.9\linewidth]{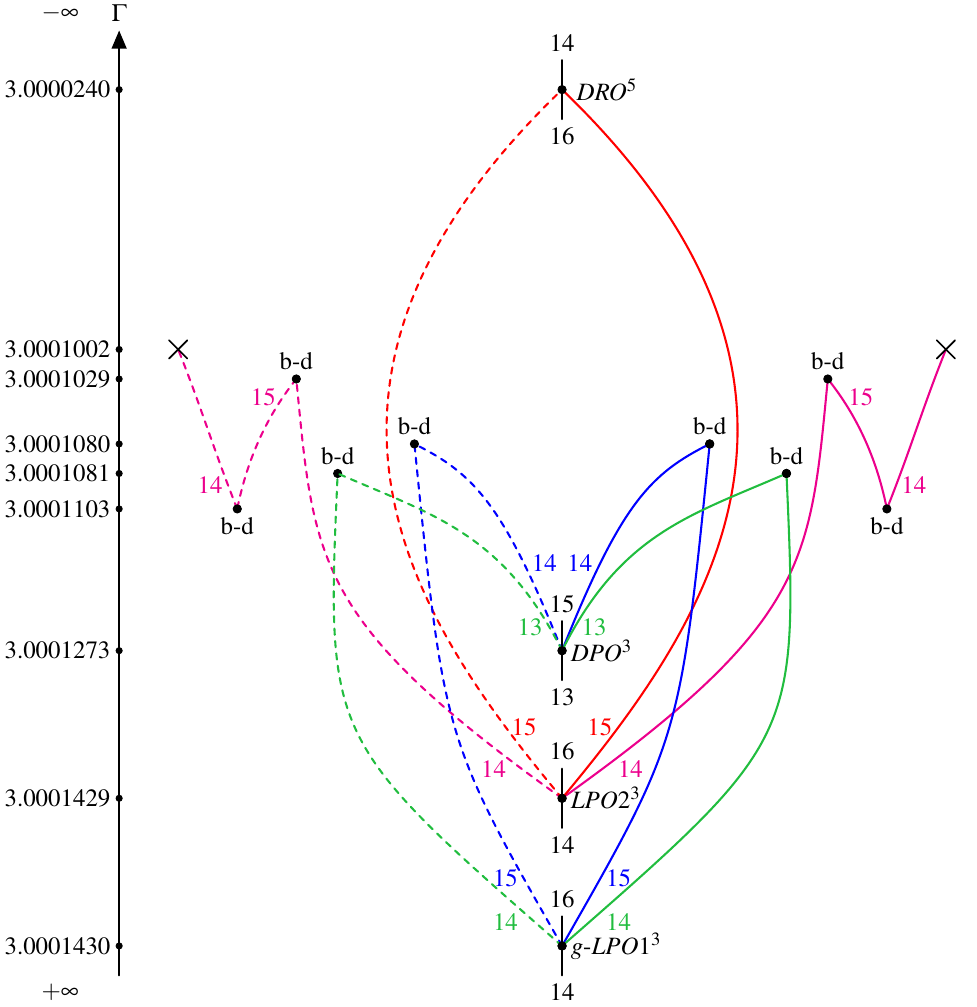}
    \caption{Bifurcation graph for Saturn-Enceladus. Its topology is exactly the same as that of Figure \ref{fig:bifurc_diagram_new}, but with different energy values.}
    \label{fig:SE_graph}
\end{figure}

\begin{figure}
    \centering
    \includegraphics[width=0.35\linewidth]{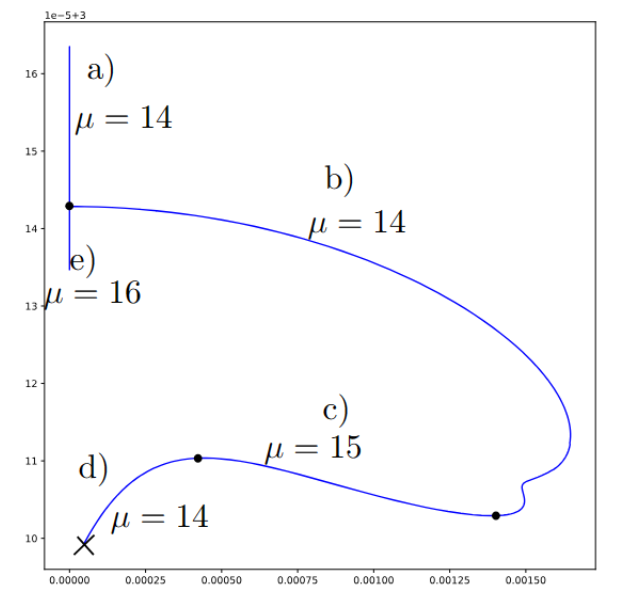}
    \caption{A bifurcation graph for Saturn-Enceladus of $xz$-plane symmetric orbits, which corresponds to the pink families in Figure \ref{fig:SE_graph}. Horizontal axis is $z$ starting value. Vertical axis is energy.}
    \label{fig:SE}
\end{figure}

\begin{figure}
        \centering
        \includegraphics[width=0.5\linewidth]{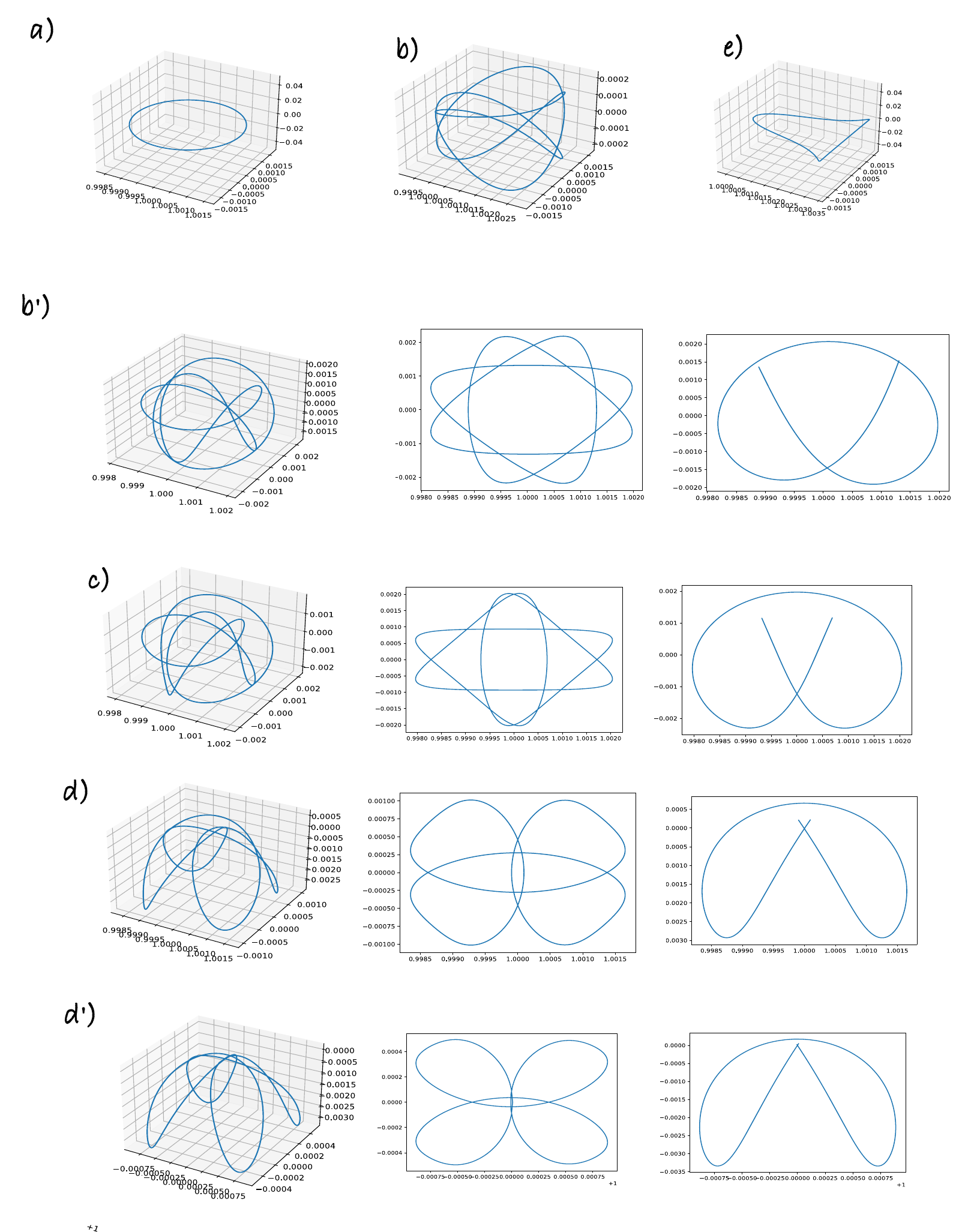}
        \caption{Plots of the orbits represented by the bifurcation graph of Figure \ref{fig:SE}. The b') orbits are part of the b) family, only for different energy: similarly for d) and d'). The second to last rows correspond to the same orbit from different angles.}
        \label{fig:SEplots}
\end{figure}

The periodic orbits in the Saturn--Enceladus system were found by continuation in the $\mu$-parameter, and its bifurcation graph corresponding to the one shown in Figure \ref{fig:bifurc_diagram_new} has exactly the same topology (but different energy values); it is plotted in Figure \ref{fig:SE_graph}. Figure \ref{fig:SE} gives a bifurcation graph corresponding to the pink families of Figure \ref{fig:SE_graph} (but drawn upside down). Note that it is not only topological, as we also record the starting value along the $z$ axis. The corresponding families of orbits are plotted in Figure \ref{fig:SEplots}.

\section{Conclusion} We presented a toolkit extracted from the general methods of symplectic geometry, aimed at studying periodic orbits of Hamiltonian systems, with their bifurcations in families, eigenvalue configurations, and stability, in a visual, and resource-efficient way. In the presence of symmetry, the information attached to orbits, and the methods involved, may be significantly refined. We illustrated these methods on numerical examples, for systems of current interest which are modelled by a restricted three-body problem (Jupiter--Europa, Saturn--Enceladus). We studied families of planar to spatial bifurcations, via bifurcation analysis and deformation from the lunar problem. The numerical findings are in agreement with the theoretical predictions, and the bifurcation graphs are completely novel. Appendix \ref{app:CZ_index_numerical} yields a self-contained documentation for the numerical implementation for computing the CZ-indices, available in GitHub. Appendix \ref{app:halo} yields an orbit in the Saturn-Enceladus system which approaches the plumes at an altitude of 29 km, and therefore may be used for future missions.

\appendix

\clearpage

\section{Numerical implementation of CZ-indices}\label{app:CZ_index_numerical}

In this appendix, we give a definition of the CZ-index due to Conley and Zehnder, needed to understand its numerical implementation. This is carried out in a Jupyter notebook, which can be found in:

\begin{verbatim}
https://github.com/ovkoert/cz-index
\end{verbatim}

We need a couple of steps to define this concept, so let's start with the goals and some intuition before getting down to the definition and computations; the CZ-index of Hamiltonian orbit is a kind of a winding number of the linearized flow along that orbit. 
\begin{figure}[htp]
\def\svgwidth{0.5\textwidth}%
\begingroup\endlinechar=-1
\resizebox{0.5\textwidth}{!}{%
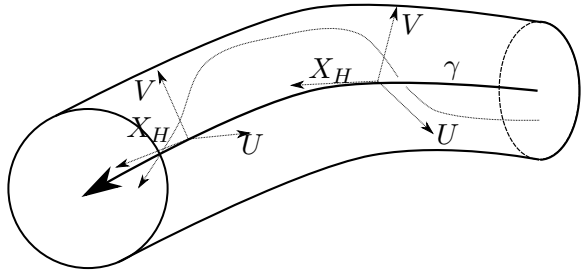%
}\endgroup
\caption{Winding of the linearized flow along an orbit with respect to a frame.}
\label{fig:linearizations_vector_fields}
\end{figure}

To remove the vagueness in this definition we need some linear algebra and a little topology. 
Consider a path of symplectic matrices, say $\psi:[0,T] \to Sp(2n)$, where $[0,T]$ is an interval of length $T$, and $Sp(2n)$ denotes the symplectic group, i.e.\
$$
Sp(2n)= \left\{ A \in M_{2n \times 2n}(\mathbb R)~|~A^t 
I  A = I \right\},
$$
with 
$$
I=\left(\begin{array}{cc}
   0  & \mathds 1 \\
   -\mathds 1  & 0
\end{array}\right).
$$
These are the matrices preserving the symplectic form $\omega(v,w)=v^tIw.$ Assume that this path starts at the identity and is non-degenerate, meaning that the endpoint $\psi(T)$ has no (generalized) eigenvalues equal to $1$.
We want to define the CZ-index $\mu_{CZ}$ as the weighted number of times the path $\psi$ goes through the eigenvalue $1$. While this can be done and makes the relation with bifurcations clearer, we choose an equivalent approach which is computationally simpler to implement. First define the Maslov cycle as the set of symplectic matrices with eigenvalue $1$, so
$$
V = \{ A \in Sp(2n)~|~ \det(A-\mathrm{id})=0 \} .
$$
The Maslov cycle $V$ divides $Sp(2n)$ into two components, namely
$$
C_+ = \{ A \in Sp(2n) ~|~\det(A-\mathrm{id})>0 \}
$$
and
$$
C_- = \{ A \in Sp(2n) ~|~\det(A-\mathrm{id})<0 \}.
$$
We choose the base points
$$
B_+ = \text{diag}(-1,-1,\ldots, -1,-1) \in C_+, \quad\text{and}\quad
B_- = \text{diag}(2,1/2,-1,-1,\ldots,-1,-1) \in C_-.
$$
We also know from the polar decomposition that any symplectic matrix $A$ can be written as $A=U S$, where $U$ is unitary and $S$ is a symmetric, positive definite matrix.
The unitary part can be extracted using the retract $\rho:Sp(2n) \to U(n) \subset Sp(2n)$,
$$
\rho(A) = (A A^T)^{-1/2}A.
$$
We can write
$$
\rho(A)=
\left(
\begin{array}{cc}
X & - Y \\
Y & X
\end{array}
\right)
,
$$
so $X+iY$ is a standard $U(n)$ matrix.
Observe that $\rho(B_-)=\mathrm{id}$ and $\rho(B_+)=-\mathrm{id}$.
Extend $\psi$ to a path $\tilde \psi:[0,2]\to Sp(2n)$ such that
\begin{itemize}
\item $\tilde \psi|_{[1,2]}$ does not intersect the Maslov cycle; and
\item $\tilde \psi(2) \in \{ B_+ ,B_- \}$.
\end{itemize} 
In other words, simply connect $\psi(1)$ to $B_+$ if $\psi(1)\in C_+$ with a path in $C_+$. Similarly if $\psi(1)\in C_-$.

Hence we get a path in the circle $U(1)$ by considering the map
$$
\gamma: t \mapsto \det{}_{\mathbb{C}}(\rho \circ \tilde \psi(t) ).
$$
This is not always a loop as it can end in $1$ or $-1$, but it will be if we double its speed. That gives us the CZ-index by taking the degree of this loop, so
$$
\mu_{CZ}(\psi) = \deg( \det{}_{\mathbb{C}}(\rho \circ \tilde \psi(t) )^2 ).
$$
Here, recall that, intuitively, the degree of a map $\gamma:[0,T]\rightarrow S^1$ taking values in the circle is the number of times it winds around the circle. So the CZ-index as defined above basically counts the number of \emph{half}-turns of the map $\gamma$ around the circle. This definition of the CZ-index can be found in \cite{Salamon1992}.

\medskip

\textbf{Extension.} The implementation simply computes the total angle change of the extension. This is straightforward to do once the right extension $\tilde \psi$ has been found. Although not very difficult, finding the extension takes up most of the script. It is based on the following theorem and observations.
\begin{theorem}
The characteristic polynomial of a symplectic matrix $A \in Sp(2n)$ is palindromic, i.e.~there are $a_0,\ldots,a_n$ such that
$$
\Delta(\lambda)=\lambda^{2n}\Delta\left(\frac{1}{\lambda}\right)=\lambda^n\sum_{k=0}^n a_k(\lambda^k+\lambda^{-k}).
$$
Furthermore, if $\lambda$ is an eigenvalue of $A$, then so are $\lambda^{-1}$, $\bar \lambda$ and $\bar \lambda^{-1}$.
\end{theorem}
This means that the eigenvalues of a non-degenerate symplectic matrix come in the following types:
\begin{itemize}
\item $(\mathcal{E})$ A pair of complex conjugate eigenvalues on the unit circle, i.e.\ an \emph{elliptic} pair;
\item ($\mathcal{H}^+$) A pair of positive real eigenvalues $a, 1/a$, i.e.\ a \emph{positive hyperbolic} pair;
\item ($\mathcal{H}^-$) A pair of negative real eigenvalues $a, 1/a$, i.e.\ a \emph{negative hyperbolic} pair;
\item ($\mathcal{C}$) A tuple of four complex eigenvalues $\lambda, 1/\lambda,\bar \lambda, 1/ \bar \lambda$ that are not real and do not lie on the unit circle, i.e\ a \emph{complex quadruple};
\end{itemize}
For the planar CR3BP, a path corresponding to the reduced monodromy will consist of $Sp(2)$-matrices. In case the path is non-degenerate, then the endpoint will be one of the following types:
\begin{itemize}
\item ($\mathcal{E}$) or ($\mathcal{H}^-$): we connect to $B_+$ and the index is odd;
\item ($\mathcal{H}^+$): we connect to $B_-$ and the index is even.
\end{itemize}
In the spatial CR3BP, a path corresponding to the reduced monodromy will be in $Sp(4)$. The following cases occur for the endpoint:
\begin{enumerate}
\item[(A)] ($\mathcal{E}^2$), ($\mathcal{EH}^-$), ($\mathcal H^{--}$), ($\mathcal{H}^{++}$), ($\mathcal{C}$): we connect to $B_+$ and the index is even;
\item[(B)] ($\mathcal{EH^+}$), ($\mathcal{H}^{-+}$): we connect to $B_-$ and the index is odd.
\end{enumerate}
Let us explain how to obtain the extension $\tilde \psi$ in the spatial case:
\begin{enumerate}
\item Eliminate all elliptic pairs;
\item Eliminate all negative hyperbolic pairs;
\item If the endpoint $S$ is of ($\mathcal{H}^{++}$) type, then take an eigenvalue decomposition $S = BDB^{-1}$, where $D$ is diagonal. The eigenvalues are generically in the form $\lambda, \mu, 1/\lambda, 1/\mu$, where $\lambda\neq \mu$. We can rescale the columns of $B$ such that they become orthogonal with respect to $\omega$. Indeed, if $S v =\lambda v$ and $Sw = \mu w$, with $v,w$ the columns of $B$, then
$$
\lambda \mu v^t I w=v^tS^tIS w =v^tI w,
$$
so $\omega(v,w)=0$, since $\lambda\mu\neq 1$. With this in mind:
\begin{itemize}
\item First deform the eigenvalues to the form $2,2,1/2,1/2$ by interpolating $D$ to $D'$;
\item Then rotate $BD'B^{-1}$ to $BR D'B^{-1}$ via a rotation matrix $R$, so that this new form is of ($\mathcal{C}$) type.
\end{itemize}
\item If the type is now ($\mathcal{C}$), then connect to $B_+$ by interpolating the eigenvalues. This finishes case (A).
\item If we are in case (B), then we have a matrix of type ($\mathcal{H}^+$) after the previous steps. By the previous observation, we may write
$$
S =B D B^{-1},
$$ 
where $B$ is a symplectic matrix and $D$ is diagonal. With the Iwasawa decomposition, we can write $B=KAN$, where $K\in U(n)$ is unitary, and $A,N$ have the form
$$
A=\left(\begin{array}{cc}
    D' & 0 \\
    0 & D'^{-1}
\end{array}\right)
$$
with $D'\in M_{2\times 2}(\mathbb R)$ diagonal and positive, and
$$
N=\left(\begin{array}{cc}
    N & M \\
    0 & N^{-t}
\end{array}\right),
$$
with $N$ upper triangular with diagonal elements equal to $1$, $NM^t=MN^t$. The matrices $K$, $A$ and $N$ can then be interpolated to the identity. The paper of Benzi and Razouk, \cite{Benzi2007}, contains an efficient and simple to implement algorithm, which we have used. 

\end{enumerate}

\medskip

\textbf{Trivializations.} We now need to connect the above linear algebra story to Hamiltonian dynamics.
Suppose that $H$ is a time-independent Hamiltonian defined on a phase space, say $M = \R^{2n}$, and consider a periodic orbit $\gamma$ of the Hamiltonian vector field $X_H$. 
We need to choose ``yard sticks'' with respect to which we measure the rotation of the linearized flow of $X_H$ as sketched in Figure \ref{fig:linearizations_vector_fields}. This is a \emph{symplectic trivialization} or \emph{frame} along the orbit $\gamma$, which simply consists of a symplectic basis of the tangent space at each point of the orbit.

There are many trivializations possible, but for the purpose of computing the CZ-index, certain choices need to be made.
For each point in the orbit $\gamma$, we take the vectors
$$
Z = \frac{1}{\Vert \nabla H \Vert^2}\nabla H,
\quad 
X_H=I \nabla H. 
$$
At each point, these two vectors span a symplectic $2$-plane $L=\langle Z, X_H \rangle$.
The normalization is chosen to ensure that
$$
\omega(Z,X_H)=dH(Z) = \frac{1}{\Vert \nabla H \Vert^2}\Vert \nabla H \Vert^2 =1.
$$
After this, we need to choose a symplectic basis of the symplectic complement $$L^\omega=\{Y: \omega(X,Y)=0 \text{ for } X \in L\}.$$
In order to get a meaningful index, we will assume that $\gamma$ is the boundary of a disk in $M$ (a \emph{spanning} disk). Then there is, by general theory, a symplectic trivialization, which is unique up to deformation. 
In particular, we can choose a symplectic basis $U_1,V_1,\ldots,U_{n-1},V_{n-1}$ of $L^\omega$. 
\begin{example}
To be concrete, let's consider the case $M=\R^4$. Rather than finding a spanning disk and searching for some trivialization, we can define a global trivialization in this case.
Let $J$ and $K$ denote the remaining two quaternionic matrices. 
For coordinates $(q_1,p_1,q_2,p_2)$, this means that
$$
I=
\left(
\begin{array}{cccc}
0&-1&0&0\\
1&0&0&0\\
0&0&0&-1\\
0&0&1&0
\end{array}
\right)
,\quad
J=\left(
\begin{array}{cccc}
0&0&1&0\\
0&0&0&-1\\
-1&0&0&0\\
0&1&0&0
\end{array}
\right)
,\quad
K=
\left(
\begin{array}{cccc}
0&0&0&1\\
0&0&1&0\\
0&-1&0&0\\
-1&0&0&0
\end{array}
\right)
.
$$ 
Then
$$
U = \frac{1}{\Vert \nabla H \Vert^2} J\nabla H,
\quad
V = K \nabla H
$$
form a symplectic trivialization of the complement $L^\omega$. This trivialization is suitable for the \emph{planar} CRTBP and also works for convex Hamiltonians, among other situations. In the Jupyter notebook, we explain how we obtain a symplectic trivialization for the \emph{spatial} CR3BP.
\end{example}
Now that we have symplectic bases of both $L_x$ and $L_x^\omega$ at a point $x$, we define the trivialization as
\[
\begin{split}
\epsilon: \gamma\times \R^{2n} & \longrightarrow \R^{2n},\\ 
(x,v) & \longmapsto
v_1 Z(x) +v_2 X_H(x) + v_3 U_1(x) +v_4 V_1(x) +\ldots  + v_{2n-1} U_{n-1}(x) +v_{2n} V_{n-1}(x).
\end{split}
\]
Now let's see how to use this trivialization.
Let $\phi^{H}_t$ denote the time $t$ flow of the Hamiltonian vector field $X_H$ and $d_x \phi^{H}_t$ the linearized flow at $x\in M$. We obtain a symplectic matrix by the formula
$$
\psi(t)=\epsilon(\phi^{H}_t(x),\cdot) ^{-1} \circ d_x\phi^{H}_t \circ \epsilon(x,\cdot).
$$
\begin{remark}
In case the notation is unfamiliar, the function {\texttt{get\_symplectic\_frame}} in the Jupyter notebook will clarify this formula. 
\end{remark}
This matrix has the form
$$
\psi=
\left(
\begin{array}{cccc}
1 & 0 & 0 & 0 \\
\delta T & 1 & \delta U & \delta V \\
c & 0 & B & \ldots \\ 
d & 0 & \ldots & \ldots 
\end{array}
\right)
$$
The matrix $B$ is a $(2n-2)\times (2n-2)$ symplectic matrix, depending on $t$, which we will call the reduced monodromy.
We will define the transverse CZ-index of the orbit $\gamma$ as the CZ-index of the path $t\mapsto B(t)$, i.e.\ $\mu_{CZ}(\gamma,\epsilon):=\mu_{CZ}(B).$
This is then the CZ-index of the orbit $\gamma$, which depends on homotopy class of trivialization. This means that the index is invariant under continuous deformations of trivializations.

\subsubsection*{Summary}
In short, the script takes initial conditions for a periodic orbit as input, and then proceeds with the following steps:
\begin{itemize}
    \item compute a numerical approximation to the orbit and the linearized flow,
    \item project the linearized flow to a symplectic frame to obtain a path of symplectic matrices,
    \item extend this path as described above, and
    \item compute the winding number of the concatenated paths.
\end{itemize}
Mathematically, the constructed path is continuous, lies in the symplectic group and the extension part of the path doesn't intersect the Maslov cycle.
However, the discretizations of both the numerical approximation and the extension have to be fine enough for the the numerical result to be correct; the script provides criteria to check this. Cases where we found that a finer discretization was necessary, include orbits that come very close to collision, and orbits that are almost degenerate.

\section{Halo orbits and polar orbits}\label{app:halo}
\begin{figure}[htb]
        \centering
        \includegraphics[width=1.0\linewidth]{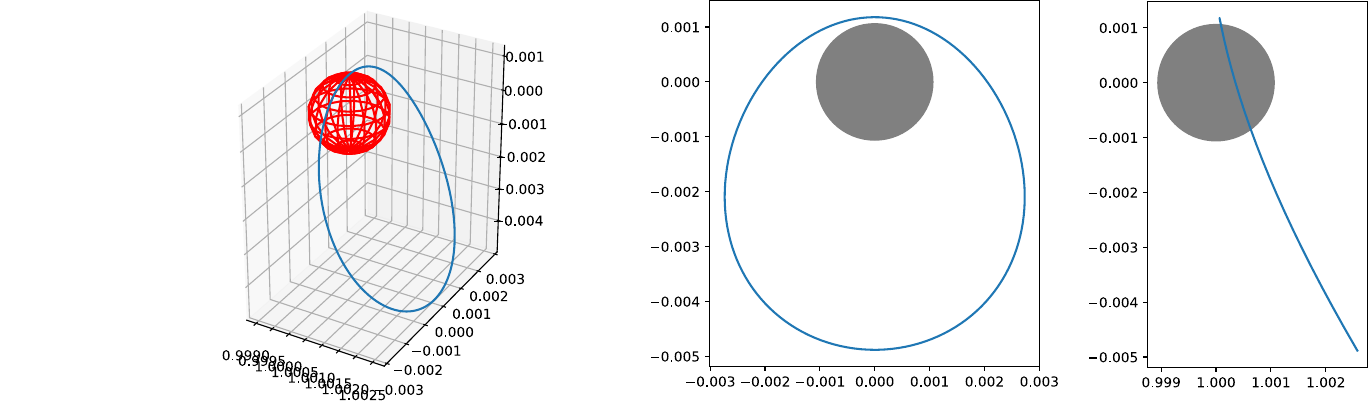}
        \caption{Plots of a Halo-polar orbit ($\Gamma=3.000034709155895$) with an altitude of 29 km. The Conley--Zehnder index has just jumped to 4, and the type is ($\mathcal{E}^2$). Europe is displayed in red and gray.}
        \label{fig:SEhalo_polar}
\end{figure}

We apply the same methods to the family of Halo orbits coming out of $L_2$ in the Saturn--Enceladus system. It turns out that this family gets very close to the plumes; the same family appears also in the whitepaper on the Enceladus Orbilander, \cite{MKGN}.
We continue this family of periodic orbits past a birth-death degeneracy to connect it to the family of polar orbits, compute the indices as well as the distance to the surface at the minimal angle with the pole.
We have taken 237,948 km for the semi-major axis of the orbit of Enceladus around Saturn, and 252.1 km for the radius of Enceladus to compute the distance to the surface. The most interesting part of the family occurs just after the index change from 3 to 4, where the orbit is both stable and getting close to the surface as illustrated in Figure~\ref{fig:SEhalo_polar}.

\begin{table}[htp]\fontsize{10}{10} \centering
\caption{Halo 2 family to polar near bifurcation for Saturn--Enceladus.}
	\scalebox{0.85}{\begin{tabular}{c|c|c|c|c|c|c}
 \hline
		$\Gamma$ & $x(0)$ & $z(0)$ & Distance to surface & Angle & $\mu_{CZ}$\\
		\hline 
        3.0000347723579006 & 1.0026629054297493 & -0.004864325835487838 & 47 km & $14^{\circ}$ & 3 \\
        3.0000347323578973 & 1.0026358028989037 & -0.004870378097113008 & 42 km & $12^{\circ}$ & 3 \\
        3.000034706717895 & 1.0025922078191964 & -0.004879135320319407 & 33 km & $10^{\circ}$ & 3 \\
        3.000034709155895 & 1.0025751548678687 & -0.004882249068671777 & 29 km & $10^{\circ}$ & 4\\
        3.000034719155896 & 1.0025570306848521 & -0.004885374862351879 & 25 km & $9^{\circ}$ & 4\\
        3.0000349743579013 & 1.0024341077005268 & -0.004902014803197094 & 0.6 km & $4^{\circ}$ & 4
        
	\end{tabular}}
	\label{halo2polar_saturn_enceladus}
\end{table}
The orbit can of course be continued further as a polar orbit, but this will result in a physical collision.

\begin{figure}
    \centering
    \includegraphics[width=0.7\linewidth]{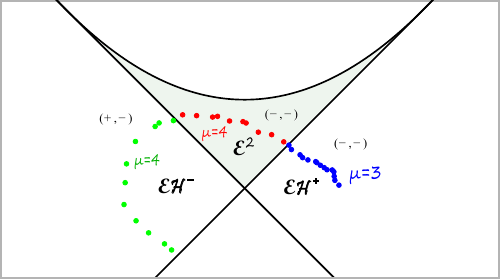}
    \caption{A GIT plot of the family of Halo orbits, with their corresponding CZ-index. The bifurcation from blue to red is of birth-death type, and from red to green, a period-doubling.}
    \label{fig:enter-label}
\end{figure}

\clearpage

\section{Tables}\label{app:tables}

In this appendix, we give tables with the data associated to the various families we have considered.

\begin{table}[htbp]\fontsize{10}{10}\selectfont \centering
\caption{Data for $g$-$LPO1$ branch for JE.}
\begin{tabular}{c|c|c|c|c|c}
\hline 
	$\Gamma$ & $x(0)$ & $\dot{y}(0)$ & $T$ & ($C/B$)-sign \& Floquet multipliers & $\mu_{CZ}^p$ / $\mu_{CZ}^s$ / $\mu_{CZ}$\\
	\hline 3.01142113 & 1.0226290 & 0.10284894 & 0.13999 & $(+/-)$ $\varphi_p = 0.137$, $(+/-)$ $\varphi_s = 0.142$ & 3 / 3 / 6\\
	3.00383366 & 1.00797270 & 0.05073828 & 1.17402 & $(+/-)$ $\varphi_p = 0.332$, $(+/-)$ $\varphi_s = 1.290$ & 3 / 3 / 6 \\
	3.00372747 & 1.00538715 & 0.07480305 & 1.33771 & $(+/-)$ $\varphi_p = 0.556$, $(+/-)$ $\varphi_s = 1.572$ & 3 / 3 / 6 \\
	3.00365597 & 1.00378076 & 0.09779235 & 1.59835 & $(+/-)$ $\varphi_p = 1.068$, $(+/-)$ $\varphi_s = 2.089$ & 3 / 3 / 6 \\
	3.00360358 & 1.00244635 & 0.12945606 & 2.08332 & $(+/-)$ $\varphi_p = 1.845$, $(+/-)$ $\varphi_s = 3.136$ & 3 / 3 / 6 \\
	3.00360326 & 1.00243628 & 0.12977936 & 2.08832 & $(+/-)$ $\varphi_p = 1.858$, $(+/+)$ $\lambda_s = -1.03$ & 3 / 3 / 6 \\
	3.00360049 & 1.00234732 & 0.13271728 & 2.13332 & $(+/-)$ $\varphi_p = 1.987$, $(+/+)$ $\lambda_s = -1.05$ & 3 / 3 / 6 \\
	3.00359960 & 1.00231829 & 0.13370950 & 2.14831 & $(+/-)$ $\varphi_p = 2.037$, $(-/+)$ $\varphi_s = 3.166$ & 3 / 3 / 6 \\
	3.00358255 & 1.00180287 & 0.15481646 & 2.43323 & $(+/+)$ $\lambda_p = -4.39$, $(-/+)$ $\varphi_s = 3.796$ & 3 / 3 / 6 \\
	3.00343430 & 1.00043030 & 0.32769866 & 3.13136 & $(+/+)$ $\lambda_p = -129$, $(-/+)$ $\varphi_s = 5.223$ & 3 / 3 / 6
\end{tabular}
\label{data_g_lpo1}
\end{table}

\begin{table}[h]\fontsize{10}{10}\selectfont \centering
\caption{Data for $DPO$ branch for JE.}
	\begin{tabular}{c|c|c|c|c|c}
 \hline
		$\Gamma$ & $x(0)$ & $\dot{y}(0)$ & $T$ & ($C/B$)-sign \& Floquet multipliers & $\mu_{CZ}^p$ / $\mu_{CZ}^s$ / $\mu_{CZ}$\\
		\hline 3.00374605 & 1.00900895 & 0.04460670 & 1.25362 & $(-/-)$ $\lambda_p = 1.29$, $(+/-)$ $\varphi_s = 1.385$ & 2 / 3 / 5\\
		3.00358658 & 1.00884026 & 0.04739922 & 1.41448 & $(-/-)$ $\lambda_p = 2.23$, $(+/-)$ $\varphi_s = 1.565$ & 2 / 3 / 5\\
		3.00356924 & 1.00889026 & 0.04727261 & 1.43426 & $(-/-)$ $\lambda_p = 2.36$, $(+/-)$ $\varphi_s = 1.589$ & 2 / 3 / 5\\
		3.00340053 & 1.00928559 & 0.04673515 & 1.64653 & $(-/-)$ $\lambda_p = 4.29$, $(+/-)$ $\varphi_s = 1.829$ & 2 / 3 / 5 \\
        3.00323697 & 1.00958786 & 0.04684116 & 1.88433 & $(-/-)$ $\lambda_p = 8.52$, $(+/-)$ $\varphi_s = 2.094$ & 2 / 3 / 5 \\
		3.00257321 & 1.00913170 & 0.05562606 & 2.88768 & $(-/-)$ $\lambda_p = 136$, $(+/-)$ $\varphi_s = 3.092$ & 2 / 3 / 5 \\
		3.00237147 & 1.00863170 & 0.05990199 & 3.16288 & $(-/-)$ $\lambda_p = 246$, $(-/+)$ $\varphi_s = 3.334$ & 2 / 3 / 5 \\
		3.00109352 & 1.00470170 & 0.09778837 & 5.12979 & $(-/-)$ $\lambda_p = 2485$, $(-/+)$ $\varphi_s = 6.161$ & 2 / 3 / 5 \\
		3.00109192 & 1.00469670 & 0.09785369 & 5.13303 & $(-/-)$ $\lambda_p = 2570$, $(+/+)$ $\lambda_s = 1.027$ & 2 / 4 / 6 \\
		3.00107109 & 1.00463170 & 0.09871030 & 5.17546 & $(-/-)$ $\lambda_p = 3062$, $(+/+)$ $\lambda_s = 1.540$ & 2 / 4 / 6 
	\end{tabular}
	\label{data_dpo}
\end{table}

\begin{table}[h]\fontsize{10}{10}\selectfont \centering
\caption{Data for $LPO2$ branch for JE.}
	\begin{tabular}{c|c|c|c|c|c}
 \hline
		$\Gamma$ & $x(0)$ & $\dot{y}(0)$ & $T$ & ($C/B$)-sign \& Floquet multipliers & $\mu_{CZ}^p$ / $\mu_{CZ}^s$ / $\mu_{CZ}$\\
		\hline 3.00374885 & 1.00955895 & 0.04118756 & 1.25694 & $(+/-)$ $\varphi_p = 0.190$, $(+/-)$ $\varphi_s = 1.397$ & 3 / 3 / 6\\
		3.00371150 & 1.01150895 & 0.03105844 & 1.34143 & $(+/-)$ $\varphi_p = 0.538$, $(+/-)$ $\varphi_s = 1.555$ & 3 / 3 / 6 \\
		3.00369790 & 1.01200895 & 0.02882949 & 1.37591 &  $(+/-)$ $\varphi_p = 0.625$, $(+/-)$ $\varphi_s = 1.620$ & 3 / 3 / 6 \\
		3.00363027 & 1.01440084 & 0.01977091 & 1.62295 & $(+/-)$ $\varphi_p = 1.101$, $(+/-)$ $\varphi_s = 2.094$ & 3 / 3 / 6 \\
		3.00357414 & 1.016776 & 0.0130372 & 2.1215 & $(+/-)$ $\varphi_p = 1.878$, $(+/-)$ $\varphi_s = 3.131$ & 3 / 3 / 6 \\
		3.00357388 & 1.016787 & 0.013014 & 2.12519 & $(+/-)$ $\varphi_p = 1.885$, $(+/+)$ $\lambda_s = -1.02$ & 3 / 3 / 6 \\
		3.00356878 & 1.01701395 & 0.01253366 & 2.20708 & $(+/-)$ $\varphi_p = 2.120$, $(-/+)$ $\varphi_s = 3.225$ & 3 / 3 / 6 \\
		3.00353952 & 1.01771395 & 0.01187914 & 2.65553 & $(+/+)$ $\lambda_p = -9.64$, $(-/+)$ $\varphi_s = 4.144$ & 3 / 3 / 6 \\
		3.00349789 & 1.01765259 & 0.01364657 & 2.95454 & $(+/+)$ $\lambda_p = -36.3$, $(-/+)$ $\varphi_s = 4.743$ & 3 / 3 / 6 
	\end{tabular}
	\label{data_lpo2}
\end{table}

\begin{table}[h]\fontsize{10}{10}\selectfont \centering
\caption{Data for $DRO$ branch for JE.}
	\begin{tabular}{c|c|c|c|c|c}
 \hline
		$\Gamma$ & $x(0)$ & $\dot{y}(0)$ & $T$ & ($C/B$)-sign \& Floquet multipliers & $\mu_{CZ}^p$ / $\mu_{CZ}^s$ / $\mu_{CZ}$\\
		\hline 3.00429783 & 0.99502455 & 0.07670173 & 0.40998 & $(-/+)$ $\varphi_p = 5.862$, $(-/+)$ $\varphi_s = 5.894$ & 1 / 1 / 2\\
		3.00156431 & 0.99037034 & 0.06224607 & 1.02778 & $(-/+)$ $\varphi_p = 5.245$, $(-/+)$ $\varphi_s = 5.406$ & 1 / 1 / 2 \\
		3.00101739 & 0.98833167 & 0.06026263 & 1.32856 & $(-/+)$ $\varphi_p = 4.973$, $(-/+)$ $\varphi_s = 5.216$ & 1 / 1 / 2 \\
		3.00060753 & 0.98623049 & 0.05949811 & 1.64998 & $(-/+)$ $\varphi_p = 4.712$, $(-/+)$ $\varphi_s = 5.050$ & 1 / 1 / 2 \\
        3.00054882 & 0.98587513 & 0.05946574 & 1.7052 & $(-/+)$ $\varphi_p = 4.670$, $(-/+)$ $\varphi_s = 5.026$ & 1 / 1 / 2 \\
		2.99962388 & 0.97762100 & 0.06369886 & 3 & $(-/+)$ $\varphi_p = 4.001$, $(-/+)$ $\varphi_s = 4.787$ & 1 / 1 / 2 \\
		2.99935885 & 0.97409965 & 0.06735824 & 3.5 & $(-/+)$ $\varphi_p = 3.987$, $(-/+)$ $\varphi_s = 4.863$ & 1 / 1 / 2 \\
		2.99908502 & 0.97038828 & 0.07212000 & 4 & $(-/+)$ $\varphi_p = 3.995$, $(-/+)$ $\varphi_s = 5.001$ & 1 / 1 / 2 \\
		2.99868251 & 0.96488658 & 0.08024713 & 4.6003 & $(-/+)$ $\varphi_p = 4.185$, $(-/+)$ $\varphi_s = 5.254$ & 1 / 1 / 2
	\end{tabular}
	
	\label{data_dro}
\end{table}

\begin{table}[h]\fontsize{10}{10} \centering
\caption{Data for one branch bifurcation from 3rd cover of the $LPO2$-orbit for JE.\ These spatial orbits are simply-symmetric w.r.t.\ the $xz$-plane and ends at collision.\ Its symmetric family is obtained by using the reflection at the $xy$-plane.}
	\scalebox{0.9}{\begin{tabular}{c|c|c|c|c|c|c}
 \hline
		$\Gamma$ & $x(0)$ & $z(0)$ & $\dot{y}(0)$ & $T$ & ($C/B$)-sign \& Floquet multipliers & $\mu_{CZ}$\\
		\hline 3.00363027 & 1.01440084 & 0 & 0.01974709 & 4.86 & $(-/+)$ $\varphi_p^3 = 3.305$, $\varphi_s^3 = 0$ & 14$\,\to\,$16 \\
        3.00362881 & 1.01439256 & 0.00046114 & 0.01976648 & 4.87 & $(-/+)$ $\varphi_1 = 3.300$, $(-/+)$ $\varphi_2 = 6.281$ & 14 \\
		3.00359018 & 1.01415816 & 0.00242577 & 0.02031752 & 4.90 & $(-/+)$ $\varphi_1 = 3.208$, $(-/+)$ $\varphi_2 = 6.280$ & 14 \\
		3.00357914 & 1.01409052 & 0.00273476 & 0.02047842 & 4.91 & $(+/+)$ $\lambda = -1.05$, $(-/+)$ $\varphi = 6.278$ & 14 \\
		3.00354287 & 1.01386628 & 0.003555363 & 0.02101794 & 4.94 & $(+/+)$ $\lambda = -1.18$, $(-/+)$ $\varphi = 6.255$ & 14 \\
		3.00325974 & 1.01198527 & 0.00688259 & 0.02594794 & 5.2 & $(+/+)$ $\lambda = -1.62$, $(-/+)$ $\varphi = 5.963$ & 14 \\
		3.00298774 & 1.00985792 & 0.00824897 & 0.03258269 & 5.5 & $(+/-)$ $\varphi_1 = 2.566$, $(-/+)$ $\varphi_2 = 5.657$ & 14 \\
		3.00270453 & 1.00652898 & 0.00795347 & 0.04651756 & 5.85 & $(-/+)$ $\varphi_1 = 1.947$, $(-/+)$ $\varphi_2 = 5.978$ & 14 \\
		3.00264234 & 1.00560524 & 0.00778449 & 0.05051319 & 5.88 & $(+/+)$ $\lambda = -1.09$, $(-/+)$ $\varphi = 5.958$ & 14 \\
		3.00263168 & 1.00544296 & 0.00774780 & 0.05124733 & 5.88 & $(-/+)$ $\varphi_1 = 3.488$, $(-/+)$ $\varphi_2 = 5.937$ & 14 \\
        3.00260038 & 1.00454296 & 0.00720347 & 0.05686831 & 5.86 & $(-/+)$ $\varphi_1 = 4.662$, $(-/+)$ $\lambda = 1$ & b-d \\
        3.00260927 & 1.00399399 & 0.00658371 & 0.06201856 & 5.8 & $(-/+)$ $\varphi = 4.813$, $(+/+)$ $\lambda = 2.278$ & 15 \\
        3.00266582 & 1.00306075 & 0.00521508 & 0.07465765 & 5.6 & $(-/+)$ $\varphi = 4.443$, $(+/+)$ $\lambda = 3.660$ & 15 \\
        3.00278841 & 1.00150186 & 0.00269606 & 0.11584022 & 5 & $(-/+)$ $\varphi = 3.653$, $(+/+)$ $\lambda = 1.829$ & 15 \\
        3.00279353 & 1.00129733 & 0.00238127 & 0.12512611 & 4.9 & $(-/+)$ $\varphi_1 = 3.540$, $(-/+)$ $\lambda = 1$ & b-d \\
        3.00277937 & 1.00084704 & 0.00169978 & 0.15351993 & 4.66 & $(-/+)$ $\varphi_1 = 3.246$, $(-/+)$ $\varphi_2 = 5.702$ & 14 
	\end{tabular}}
	\label{data_lpo2_3rd_cover_1}
\end{table}

\begin{table}[t]\fontsize{10}{10} \centering
\caption{Data for one branch bifurcation from $LPO2^3$-orbit for JE, ending at $DRO^5$.\ The CZ-index is constant, and gives a bridge between the planar orbits.\ These spatial orbits are $x$-axis-symmetric. Its symmetric family is obtained by reflection at the ecliptic.}
	\scalebox{0.9}{\begin{tabular}{c|c|c|c|c|c|c}
 \hline
		$\Gamma$ & $x(0)$ & $\dot{y}(0)$ & $\dot{z}(0)$ & $T$ & ($C/B$)-sign \& Floquet multipliers & $\mu_{CZ}$\\
		\hline 3.00363027 & 1.01440084 & 0.01974709 & 0 & 4.86 & $(-/+)$ $\varphi_p^3 = 3.305$, $\varphi_s^3 = 0$ & 14$\,\to\,$16 \\
        3.00351924 & 1.01408954 & 0.01928185 & 0.01342237 & 4.96 & $(-/-)$ $\lambda_1 = -1.24$, $(-/-)$ $\lambda_2 = 1.09$ & 15 \\
		3.00321170 & 1.01314307 & 0.01799332 & 0.02676988 & 5.25 & $(-/-)$ $\lambda_1 = -1.62$, $(-/-)$ $\lambda_2 = 1.45$ & 15 \\
		3.00302231 & 1.01246670 & 0.01723541 & 0.03300118 & 5.46 & $(+/-)$ $\varphi = 2.749$, $(-/-)$ $\lambda = 1.78$ & 15 \\
		3.00273486 & 1.01099334 & 0.01657953 & 0.04285920 & 5.82 & $(+/-)$ $\varphi = 1.618$, $(-/-)$ $\lambda = 1.55$ & 15 \\
		3.00270684 & 1.01077857 & 0.01644568 & 0.04412031 & 5.85 & $(+/-)$ $\varphi = 1.801$, $(-/-)$ $\lambda = 1.40$ & 15 \\
		3.00266563 & 1.01548137 & 0.01548137 & 0.04575934 & 5.88 & $(+/-)$ $\varphi = 2.506$, $(-/-)$ $\lambda = 1.35$ & 15 \\
		3.00243536 & 1.01068879 & 0.00516070 & 0.04995370 & 6 & $(-/+)$ $\varphi = 5.761$, $(-/-)$ $\lambda = 8.52$ & 15 \\
		3.00204821 & 1.01119864 & $-$0.01174966 & 0.05089250 & 6.24 & $(-/+)$ $\varphi = 5.965$, $(-/-)$ $\lambda = 31.1$ & 15 \\
		3.00172312 & 1.01167768 & $-$0.02457421 & 0.04793350 & 6.5 & $(-/+)$ $\varphi = 6.016$, $(-/-)$ $\lambda = 30.0$ & 15 \\
		3.00147493 & 1.01207539 & $-$0.03349482 & 0.04374826 & 6.75 & $(-/+)$ $\varphi = 6.070$, $(-/-)$ $\lambda = 22.3$ & 15 \\
		3.00127220 & 1.01242785 & $-$0.04020652 & 0.03910792 & 7 & $(-/+)$ $\varphi = 6.124$, $(-/-)$ $\lambda = 15.2$ & 15 \\
		3.00096072 & 1.01304236 & $-$0.04944170 & 0.02947258 & 7.5 & $(-/+)$ $\varphi = 6.201$, $(-/-)$ $\lambda = 6.16$ & 15 \\
		3.00073221 & 1.01358551 & $-$0.05528744 & 0.01932635 & 8 & $(-/+)$ $\varphi = 5.911$, $(-/-)$ $\lambda = 1.14$ & 15\\
		3.00055690 & 1.01409401 & $-$0.05914769 & 0.00388381 & 8.5 & $(-/+)$ $\varphi = 4.550$, $(-/-)$ $\lambda = 1.00$ & 15 \\
        3.00054882 & 1.01412064 & $-$0.05930512 & 0 & 8.52 & $(-/+)$ $\varphi_p^5 = 4.500$, $\varphi_s^5 = 0$ & 16$\,\to\,$14 \\
	\end{tabular}}
	\label{data_lpo2_3rd_cover_2}
\end{table}

\begin{table}[b]\fontsize{10}{10} \centering
\caption{Data for one branch bifurcation from $g$-$LPO1^3$ for JE.\ These spatial orbits are $x$-axis-symmetric and connected to one branch bifurcation from $DPO^3$ via b-d.\ Its symmetric family is obtained by reflection at the ecliptic.}
	\scalebox{0.9}{\begin{tabular}{c|c|c|c|c|c|c}
 \hline
		$\Gamma$ & $x(0)$ & $\dot{y}(0)$ & $\dot{z}(0)$ & $T$ & ($C/B$)-sign \& Floquet multipliers & $\mu_{CZ}$\\
		\hline 3.00365597 & 0.98557900 & $-$0.01951876 & 0 & 4.79 & $(-/+)$ $\varphi_p^3 = 3.205$, $\varphi_s^3 = 0$ & 14$\,\to\,$16 \\
        3.00365338 & 0.98556744 & $-$0.01945341 & 0.00182488 & 4.80 & $(-/+)$ $\varphi = 3.214$, $(-/-)$ $\lambda = 1.001$ & 15 \\
		3.00363389 & 0.98561874 & $-$0.01938011 & 0.00580611 & 4.82 & $(-/-)$ $\lambda_1 = -1.01$, $(-/-)$ $\lambda_2 = 1.005$ & 15 \\
        3.00329911 & 0.98657072 & $-$0.01815373 & 0.02416862 & 5.12 & $(-/-)$ $\lambda_1 = -1.54$, $(-/-)$ $\lambda_2 = 1.378$ & 15 \\
        3.00314093 & 0.98708523 & $-$0.01763006 & 0.02950246 & 5.30 & $(-/-)$ $\lambda_1 = -1.02$, $(-/-)$ $\lambda_2 = 1.685$ & 15 \\
        3.00300399 & 0.98759083 & $-$0.01727225 & 0.03377964 & 5.46 & $(+/-)$ $\varphi = 2.290$, $(-/-)$ $\lambda = 1.977$ & 15 \\
        3.00281046 & 0.98856773 & $-$0.01745284 & 0.04009001 & 5.73 & $(+/-)$ $\varphi = 0.976$, $(-/-)$ $\lambda = 2.451$ & 15 \\
        3.00275889 & 0.98918471 & $-$0.01885469 & 0.04265726 & 5.82 & $(+/-)$ $\varphi = 0.134$, $(-/-)$ $\lambda = 4.422$ & 15 \\
         & & & & & & b-d \\
         3.00275823 & 0.98925887 & $-$0.01917383 & 0.04284561 & 5.82 & $(-/-)$ $\lambda_1 = 1.041$, $(-/-)$ $\lambda_2 = 5.013$ & 14 \\
         3.00276196 & 0.98939330 & $-$0.01992825 & 0.04305515 & 5.83 & $(-/-)$ $\lambda_1 = 1.148$, $(-/-)$ $\lambda_2 = 6.636$ & 14 \\
         3.00296320 & 0.98997681 & $-$0.03157520 & 0.03598567 & 5.75 & $(-/-)$ $\lambda_1 = 1.064$, $(-/-)$ $\lambda_2 = 100.2$ & 14 \\
         3.00316033 & 0.99025736 & $-$0.04249067 & 0.02045809 & 5.68 & $(-/-)$ $\lambda_1 = 1.011$, $(-/-)$ $\lambda_2 = 400.3$ & 14 \\
         3.00323676 & 0.99035914 & $-$0.04685850 & 0.00096495 & 5.65 & $(-/-)$ $\lambda_1 = 1.0001$, $(-/-)$ $\lambda_2 = 619.3$ & 14 \\
         3.00323697 & 0.99035942 & $-$0.04686768 & 0 & 5.65 & $\varphi_s^3 = 0$, $(-/-)$ $\lambda = 620.23$ & 13$\,\to\,$15
	\end{tabular}}
	\label{data_lpo1_3rd_cover_1}
\end{table}

%\vspace{-15cm}

\begin{table}[h]\fontsize{10}{10} \centering
\caption{Data for one branch bifurcation from 3rd cover of the $g$-$LPO1$-orbit for JE.\ These spatial orbits are simply-symmetric w.r.t.\ the $xz$-plane and they are connected to one branch bifurcation from the 3rd cover of the $DPO$-orbit via birth-death.\ Its symmetric family is obtained by using the reflection at the ecliptic.}
	\scalebox{0.9}{\begin{tabular}{c|c|c|c|c|c|c}
 \hline
		$\Gamma$ & $x(0)$ & $z(0)$ & $\dot{y}(0)$ & $T$ & ($C/B$)-sign \& Floquet multipliers & $\mu_{CZ}$\\
		\hline 3.00365597 & 0.98557900 & 0 & $-$0.01951876 & 4.79 & $(-/+)$ $\varphi_p^3 = 3.205$, $\varphi_s^3 = 0$ & 14$\,\to\,$16 \\
        3.00365461 & 0.98556706 & $-$0.00036219 & $-$0.01947401 & 4.80 & $(-/+)$ $\varphi_1 = 3.215$, $(-/+)$ $\varphi_2 = 6.281$ & 14 \\
		3.00363389 & 0.98568597 & $-$0.00175392 & $-$0.01975974 & 4.82 & $(+/+)$ $\lambda = -1.01$, $(-/+)$ $\varphi = 6.277$ & 14 \\
		3.00360033 & 0.98588066 & $-$0.00278699 & $-$0.02023321 & 4.84 & $(+/+)$ $\lambda = -1.12$, $(-/+)$ $\varphi = 6.263$ & 14 \\
        3.00331461 & 0.98766211 & $-$0.00655464 & $-$0.02492170 & 5.11 & $(+/+)$ $\lambda = -1.52$, $(-/+)$ $\varphi = 5.978$ & 14 \\
        3.00314742 & 0.98883839 & $-$0.00763969 & $-$0.02842198 & 5.29 & $(+/-)$ $\varphi_1 = 3.031$, $(-/+)$ $\varphi_2 = 5.756$ & 14 \\
        3.00289637 & 0.99094696 & $-$0.00836889 & $-$0.03579903 & 5.61 & $(+/-)$ $\varphi_1 = 1.373$, $(-/+)$ $\varphi_2 = 5.309$ & 14 \\
        3.00285045 & 0.99142732 & $-$0.00835072 & $-$0.03776509 & 5.67 & $0.376 \pm 0.570 i$, $0.806 \pm 1.221 i$ & 14 \\
        3.00277633 & 0.99243798 & $-$0.00805049 & $-$0.04244268 & 5.78 & $0.491 \pm 0.121 i$, $1.917 \pm 0.472 i$ & 14 \\
        3.00277358 & 0.99249304 & $-$0.00802039 & $-$0.04284315 & 5.79 & $(+/+)$ $\lambda_1 = 1.987$, $(-/-)$ $\lambda_2 = 2.016$ & 14 \\
        3.00276770 & 0.993012244 & $-$0.00750257 & $-$0.04644237 & 5.82 & $(+/+)$ $\lambda_1 = 1.000$, $(-/-)$ $\lambda_2 = 7.181$ & 14 \\
         & & & & & & b-d \\
        3.00277093 & 0.99302677 & 0.00744432 & $-$0.04668427 & 5.82 & $(-/+)$ $\varphi = 6.138$, $(-/-)$ $\lambda = 8.013$ & 13 \\
        3.00296373 & 0.99201168 & 0.00567451 & $-$0.04755480 & 5.75	& $(-/+)$ $\varphi = 6.233$, $(-/-)$ $\lambda = 100.3$ & 13 \\
        3.00316033 & 0.99081340 & 0.00306250 & $-$0.04704925 & 5.68 & $(-/+)$ $\varphi = 6.263$, $(-/-)$ $\lambda = 400.3$ & 13 \\
        3.00323697 & 0.99035942 & 0 & $-$0.04686768 & 5.65 & $(-/-)$ $\lambda = 620.23$, $\varphi_s^3 = 0$ & 13$\,\to\,$15
	\end{tabular}}
	\label{data_lpo1_3rd_cover_2}
\end{table}

\begin{table}[htbp]\fontsize{10}{10}\selectfont \centering
\caption{$g$-$LPO1$ branch for SE.}
\begin{tabular}{c|c|c|c|c|c}
\hline 
	$\Gamma$ & $x(0)$ & $\dot{y}(0)$ & $T$ & ($C/B$)-sign \& Floquet multipliers & $\mu_{CZ}^p$ / $\mu_{CZ}^s$ / $\mu_{CZ}$\\
	\hline 3.00033109 & 1.00061213 & 0.01702742 & 0.22681 & $(+/-)$ $\varphi_p = 0.218$, $(+/-)$ $\varphi_s = 0.232$ & 3 / 3 / 6\\
    3.00015209 & 1.00157026 & 0.00982528 & 1.13552 & $(+/-)$ $\varphi_p = 0.421$, $(+/-)$ $\varphi_s = 1.245$ & 3 / 3 / 6\\
    3.00014609 & 1.00109981 & 0.01422219 & 1.33077 & $(+/-)$ $\varphi_p = 0.534$, $(+/-)$ $\varphi_s = 1.549$ & 3 / 3 / 6 \\
    3.00014309 & 1.00076095 & 0.01889960 & 1.60328 & $(+/-)$ $\varphi_p = 1.118$, $(+/-)$ $\varphi_s = 2.081$ & 3 / 3 / 6 \\
    3.00014089 & 1.00047925 & 0.02552883 & 2.13798 & $(+/-)$ $\varphi_p = 1.962$, $(+/+)$ $\lambda_s = -1.06$ & 3 / 3 / 6 \\
    3.00014069 & 1.00044659 & 0.02664172 & 2.22579 & $(+/-)$ $\varphi_p = 2.302$, $(-/+)$ $\varphi_s = 3.324$ & 3 / 3 / 6 \\
    3.00014049 & 1.00041446 & 0.02785333 & 2.31612 & $(+/+)$ $\lambda_p = -1.58$, $(-/+)$ $\varphi_s = 3.513$ & 3 / 3 / 6 \\
    3.00013817 & 1.00020619 & 0.04126365 & 2.90065 & $(+/+)$ $\lambda_p = -35.5$, $(-/+)$ $\varphi_s = 4.711$ & 3 / 3 / 6
\end{tabular}
\label{data_g_lpo1_s_e}
\end{table}

\begin{table}[h]\fontsize{10}{10}\selectfont \centering
\caption{$DPO$ branch for SE.}
	\begin{tabular}{c|c|c|c|c|c}
 \hline
		$\Gamma$ & $x(0)$ & $\dot{y}(0)$ & $T$ & ($C/B$)-sign \& Floquet multipliers & $\mu_{CZ}^p$ / $\mu_{CZ}^s$ / $\mu_{CZ}$\\
		\hline 3.00014744 &	0.99838904 & $-$0.00977411 & 1.23860 & $(-/-)$ $\lambda_p = 1.188$, $(+/-)$ $\varphi_s = 1.364$ & 2 / 3 / 5\\
        3.00014064 & 0.99826971 & $-$0.00934116 & 1.41906 & $(-/-)$ $\lambda_p = 2.270$, $(+/-)$ $\varphi_s = 1.571$ & 2 / 3 / 5\\
        3.00012744 & 0.99811661 & $-$0.00917943 & 1.88166 & $(-/-)$ $\lambda_p = 8.463$, $(+/-)$ $\varphi_s = 2.091$ & 2 / 3 / 5\\
        3.00011304 & 0.99809959 & $-$0.00985129 & 2.46619 & $(-/-)$ $\lambda_p = 46.25$, $(+/-)$ $\varphi_s = 2.697$ & 2 / 3 / 5\\
        3.00010524 & 0.99815786 & $-$0.01051283 & 2.76504 & $(-/-)$ $\lambda_p = 101.1$, $(+/-)$ $\varphi_s = 2.979$ & 2 / 3 / 5\\
        3.00008192 & 0.99846180 & $-$0.01309687 & 3.59018 & $(-/-)$ $\lambda_p = 500.3$, $(+/-)$ $\varphi_s = 3.699$ & 2 / 3 / 5\\
        3.00006838 & 0.99867270 & $-$0.01491712 & 4.08764 & $(-/-)$ $\lambda_p = 911.4$, $(+/-)$ $\varphi_s = 4.164$ & 2 / 3 / 5
        \end{tabular}
	\label{data_dpo_s_e}
\end{table}

\begin{table}[h]\fontsize{10}{10}\selectfont \centering
\caption{$LPO2$ branch for SE.}
	\begin{tabular}{c|c|c|c|c|c}
 \hline
		$\Gamma$ & $x(0)$ & $\dot{y}(0)$ & $T$ & ($C/B$)-sign \& Floquet multipliers & $\mu_{CZ}^p$ / $\mu_{CZ}^s$ / $\mu_{CZ}$\\
		\hline 3.00014639 & 1.00217453 & 0.00646333 & 1.30672 & $(+/-)$ $\varphi_p = 0.451$, $(+/-)$ $\varphi_s = 1.499$ & 3 / 3 / 6\\
        3.00014319 & 1.00276240 & 0.00405869 & 1.57114 & $(+/-)$ $\varphi_p = 1.018$, $(+/-)$ $\varphi_s = 2.009$ & 3 / 3 / 6\\
        3.00014299 & 1.00279991 & 0.00393120 & 1.59627 & $(+/-)$ $\varphi_p = 1.060$, $(+/-)$ $\varphi_s = 2.059$ & 3 / 3 / 6\\
        3.00014075 & 1.00328757 & 0.00253817 & 2.11312 & $(+/-)$ $\varphi_p = 1.881$, $(+/+)$ $\lambda_s = -1.03$ & 3 / 3 / 6\\
        3.00014061 & 1.00331899 & 0.00247078 & 2.17031 & $(+/-)$ $\varphi_p = 2.043$, $(-/+)$ $\varphi_s = 3.142$ & 3 / 3 / 6 \\
		  & & & & & \\
		3.00014030 & 1.00338269 & 0.00235356 & 2.31077 & $(+/+)$ $\lambda_p = -1.22$, $(-/+)$ $\varphi_s = 3.484$ & 3 / 3 / 6 \\
		3.00013613 & 1.00339866 & 0.00308260 & 3.07749 & $(+/+)$ $\lambda_p = -69.2$, $(-/+)$ $\varphi_s = 5.041$ & 3 / 3 / 6
	\end{tabular}
	\label{data_lpo2_s_e}
\end{table}

\begin{table}[h]\fontsize{10}{10}\selectfont \centering
\caption{$DRO$ branch for SE.}
	\begin{tabular}{c|c|c|c|c|c}
 \hline
		$\Gamma$ & $x(0)$ & $\dot{y}(0)$ & $T$ & ($C/B$)-sign \& Floquet multipliers & $\mu_{CZ}^p$ / $\mu_{CZ}^s$ / $\mu_{CZ}$\\
		\hline 3.00010525 & 1.00137692 & $-$0.01325298 & 0.67023 & $(-/+)$ $\varphi_p = 5.595$, $(-/+)$ $\varphi_s = 5.673$ & 1 / 1 / 2\\
        3.00004405 & 1.00224224 & $-$0.01182054 & 1.29643 & $(-/+)$ $\varphi_p = 5.001$, $(-/+)$ $\varphi_s = 5.234$ & 1 / 1 / 2 \\
        3.00002425 & 1.00276541 & $-$0.01163206 & 1.70339 & $(-/+)$ $\varphi_p = 4.671$, $(-/+)$ $\varphi_s = 5.027$ & 1 / 1 / 2 \\
        2.99999205 & 1.00419386 & $-$0.01226462 & 2.84090 & $(-/+)$ $\varphi_p = 4.060$, $(-/+)$ $\varphi_s = 4.785$ & 1 / 1 / 2 \\
        2.99996645 & 1.00592895 & $-$0.01421020 & 4.06771 & $(-/+)$ $\varphi_p = 4.007$, $(-/+)$ $\varphi_s = 5.025$ & 1 / 1 / 2 \\
        2.99995365 & 1.00682933 & $-$0.01550883 & 4.56302 & $(-/+)$ $\varphi_p = 4.172$, $(-/+)$ $\varphi_s = 5.236$ & 1 / 1 / 2 \\
        2.99986545 & 1.01156940 & $-$0.02377257 & 5.79878 & $(-/+)$ $\varphi_p = 5.114$, $(-/+)$ $\varphi_s = 5.951$ & 1 / 1 / 2
	\end{tabular}
	\label{data_dro_s_e}
\end{table}

\begin{table}[h]\fontsize{10}{10} \centering
\caption{Purple 14 for SE. These spatial orbits are simply-symmetric w.r.t.\ the $xz$-plane and ends at collision.\ Its symmetric family is obtained by using the reflection at the $xy$-plane.}
	\scalebox{0.9}{\begin{tabular}{c|c|c|c|c|c|c}
 \hline
		$\Gamma$ & $x(0)$ & $z(0)$ & $\dot{y}(0)$ & $T$ & ($C/B$)-sign \& Floquet multipliers & $\mu_{CZ}$\\
		\hline 3.00014299 & 1.00279991 & 0 & 0.00393120 & 4.78 & $(-/+)$ $\varphi_p^3 = 3.182$, $\varphi_s^3 = 0$ & 14$\,\to\,$16 \\
        3.00014260 & 1.00281790 & 0.00019127 & 0.00386379 & 4.84 & $(-/+)$ $\varphi_1 = 3.259$, $(-/+)$ $\varphi_2 = 6.282$ & 14 \\
        3.00014129 & 1.00277804 & 0.00047538 & 0.00395858 & 4.87 & $(+/+)$ $\lambda = 3.208$, $(-/+)$ $\varphi = 6.280$ & 14 \\
        3.00012209 & 1.00211470 & 0.00153818 & 0.00578987 & 5.36 & $(+/-)$ $\varphi_1 = 3.111$, $(-/+)$ $\varphi_2 = 5.740$ & 14 \\
        3.00010829 & 1.00137563 & 0.00156298 & 0.00875806 &	5.82 & $(-/+)$ $\varphi_1 = 0.776$, $(-/+)$ $\varphi_2 = 5.669$ & 14 \\
        3.00010291 & 1.00087722 & 0.00139401 & 0.01128814 & 5.86 & $(+/-)$ $\varphi_1 = 4.665$, $(-/+)$ $\varphi_2 = 6.126$ & 14 \\
        & & & & & & b-d \\
        3.00010295 & 1.00083899 & 0.00135458 & 0.01162196 & 5.84 & $(-/+)$ $\varphi = 4.783$, $(+/+)$ $\lambda = 1.517$ & 15 \\
        3.00010612 & 1.00055789 & 0.00095022 & 0.01542803 & 5.53 & $(-/+)$ $\varphi = 4.334$, $(+/+)$ $\lambda = 3.712$ & 15 \\
        3.00011036 & 1.00024558 & 0.00044664 & 0.02518603 & 4.86 & $(-/+)$ $\varphi = 3.504$, $(+/+)$ $\lambda = 1.109$ & 15 \\
        & & & & & & b-d \\
        3.00011031 & 1.00021998 & 0.00040764 & 0.02663091 & 4.80 & $(-/+)$ $\varphi_1 = 3.424$, $(-/+)$ $\varphi_2 = 5.916$ & 14 \\
        3.00010281 & 1.00003238 & 0.00010169 & 0.05878312 & 4.17 & $(-/+)$ $\varphi_1 = 3.781$, $(-/+)$ $\varphi_2 = 5.884$ & 14 \\
        3.00010021 & 1.00001491 & 0.00006054 & 0.07739259 & 4.08 & $(-/+)$ $\varphi_1 = 3.965$, $(-/+)$ $\varphi_2 = 6.041$ & 14
	\end{tabular}}
	\label{purple14_saturn_enceladus}
\end{table}

\begin{table}[h]\fontsize{10}{10} \centering
\caption{Red 15 for SE. The CZ-index is constant, and gives a bridge between the planar orbits.\ These spatial orbits are $x$-axis-symmetric. Its symmetric family is obtained by reflection at the ecliptic.}
	\scalebox{0.9}{\begin{tabular}{c|c|c|c|c|c|c}
 \hline
		$\Gamma$ & $x(0)$ & $\dot{y}(0)$ & $\dot{z}(0)$ & $T$ & ($C/B$)-sign \& Floquet multipliers & $\mu_{CZ}$\\
		\hline 3.00014299 & 1.00279991 & 0.00393120 & 0 & 4.78 & $(-/+)$ $\varphi_p^3 = 3.182$, $\varphi_s^3 = 0$ & 14$\,\to\,$16 \\
        3.00013559 & 1.00272159 & 0.00369764 & 0.00343416 & 5.00 & $(-/-)$ $\lambda_1 = -1.39$, $(-/-)$ $\lambda_2 = 1.13$ & 15 \\
        3.00012559 & 1.00256280 & 0.00349797 & 0.00541943 & 5.26 & $(-/-)$ $\lambda_1 = -1.53$, $(-/-)$ $\lambda_2 = 1.52$ & 15 \\
        3.00012139 & 1.00248791 & 0.00342004 & 0.00611889 & 5.38 & $(+/-)$ $\varphi = 2.928$, $(-/-)$ $\lambda = 1.73$ & 15 \\
        3.00011759 & 1.00241319 & 0.00335832 & 0.00672992 & 5.50 & $(+/-)$ $\varphi = 2.395$, $(-/-)$ $\lambda = 1.90$ & 15 \\
        3.00010399 & 1.00204385 & 0.00282460 & 0.00925842 & 5.90 & $(-/+)$ $\varphi = 3.184$, $(-/-)$ $\lambda = 1.41$ & 15 \\
        3.00009759 & 1.00208424 & 0.00127533 & 0.00976558 & 5.98 & $(-/+)$ $\varphi = 5.692$, $(-/-)$ $\lambda = 6.55$ & 15 \\
        3.00009417 & 1.00210781 & 0.00047162 & 0.00992209 & 6.03 & $(-/+)$ $\varphi = 5.846$, $(-/-)$ $\lambda = 13.5$ & 15 \\
        3.00009145 & 1.00212691 & $-$0.00015208 & 0.00999912 & 6.07 & $(-/+)$ $\varphi = 5.888$, $(-/-)$ $\lambda = 18.8$ & 15 \\
        3.00007425 & 1.00225503 & $-$0.00382600 & 0.00970041 & 6.38 & $(-/+)$ $\varphi = 5.974$, $(-/-)$ $\lambda = 32.3$ & 15 \\
        3.00005185 & 1.00244477 & $-$0.00788343 & 0.00765368 & 7.00 & $(-/+)$ $\varphi = 6.092$, $(-/-)$ $\lambda = 15.2$ & 15 \\
        3.00004205 & 1.00254154 & $-$0.00938862 & 0.00615688 & 7.40 & $(-/+)$ $\varphi = 6.171$, $(-/-)$ $\lambda = 7.51$ & 15 \\
        3.00003105 & 1.00266948 & $-$0.01085595 & 0.00374582 & 8.01 & $(-/+)$ $\varphi = 5.829$, $(-/-)$ $\lambda = 1.08$ & 15 \\
        3.00002425 & 1.00276636 & $-$0.01162574 & 0.00034005 & 8.51 & $(-/+)$ $\varphi = 4.510$, $(-/-)$ $\lambda = 1.05$ & 15\\
        3.00002405 & 1.00277176 & $-$0.01163165 & 0 &	8.54 & $(-/+)$ $\varphi_p^5 = 4.489$, $\varphi_s^5 = 0$ & 16$\,\to\,$14 \\  
	\end{tabular}}
	\label{red15_saturn_enceladus}
\end{table}

\begin{table}[h]\fontsize{10}{10} \centering
\caption{Blue 15 and Blue 14 for SE. These spatial orbits are $x$-axis-symmetric and connected to one branch bifurcation from $DPO^3$ via b-d.\ Its symmetric family is obtained by reflection at the ecliptic.}
	\scalebox{0.9}{\begin{tabular}{c|c|c|c|c|c|c}
 \hline
		$\Gamma$ & $x(0)$ & $\dot{y}(0)$ & $\dot{z}(0)$ & $T$ & ($C/B$)-sign \& Floquet multipliers & $\mu_{CZ}$\\
		\hline 3.00014309 & 0.99718487 & $-$0.00386895 & 0 & 4.80 & $(-/+)$ $\varphi_p^3 = 3.217$, $\varphi_s^3 = 0$ & 14$\,\to\,$16 \\
        3.00014289 & 0.99717670 & $-$0.00383214 & 0.00048106 & 4.83 & $(-/+)$ $\varphi = 3.246$, $(-/-)$ $\lambda = 1.016$ & 15 \\
        3.00014169 & 0.99719319 & $-$0.00380810 & 0.00145833 & 4.85 & $(-/-)$ $\lambda_1 = -1.04$, $(-/-)$ $\lambda_2 = 1.027$ & 15 \\
        3.00013425 & 0.99730021 & $-$0.00365996 & 0.00377737 & 5.02 & $(-/-)$ $\lambda_1 = -1.45$, $(-/-)$ $\lambda_2 = 1.183$ & 15 \\
        3.00012251 & 0.99749232 & $-$0.00343875 & 0.00594605 & 5.34 & $(-/-)$ $\lambda_1 = -1.03$, $(-/-)$ $\lambda_2 = 1.700$ & 15 \\
        3.00012093 & 0.99752153 & $-$0.00341248 & 0.00620193 & 5.39 &	$(+/-)$ $\varphi = 2.758$, $(-/-)$ $\lambda = 1.782$ & 15 \\
        3.00011193 & 0.99772535 & $-$0.00333595 & 0.00766874 & 5.69 & $(+/-)$ $\varphi = 0.976$, $(-/-)$ $\lambda = 2.451$ & 15 \\
        3.00010800 & 0.99793026 & $-$0.00365343 & 0.00862538 & 5.84 & $(+/-)$ $\varphi = 0.092$, $(-/-)$ $\lambda = 2.549$ & 15 \\
         & & & & & & b-d \\
        3.00010820 & 0.99796166 & $-$0.00381731 & 0.00868408 & 5.84 & $(-/-)$ $\lambda_1 = 1.196$, $(-/-)$ $\lambda_2 = 3.812$ & 14 \\
        3.00010958 & 0.99799080 & $-$0.00423992 & 0.00854277 & 5.83 & $(-/-)$ $\lambda_1 = 1.171$, $(-/-)$ $\lambda_2 = 10.02$ & 14 \\ 
        3.00012042 & 0.99807079 & $-$0.00718764 & 0.00594564 & 5.71 & $(-/-)$ $\lambda_1 = 1.068$, $(-/-)$ $\lambda_2 = 203.9$ & 14 \\ 
        3.00012568 & 0.99810519 & $-$0.00868716 & 0.00307711 & 5.66 & $(-/-)$ $\lambda_1 = 1.011$, $(-/-)$ $\lambda_2 = 485.0$ & 14 \\ 
        3.00012738 & 0.99811615 & $-$0.00918031 & 0 & 5.65 & $\varphi_s^3 = 0$, $(-/-)$ $\lambda_p = 620.1$ & 13$\,\to\,$15
	\end{tabular}}
	\label{blue15_saturn_enceladus}
\end{table}

\begin{table}[h]\fontsize{10}{10} \centering
\caption{Green 14 and Green 13 for SE. These spatial orbits are simply-symmetric w.r.t.\ the $xz$-plane and they are connected to one branch bifurcation from the 3rd cover of the $DPO$-orbit via b-d.\ Its symmetric family is obtained by using the reflection at the ecliptic.}
	\scalebox{0.9}{\begin{tabular}{c|c|c|c|c|c|c}
 \hline
		$\Gamma$ & $x(0)$ & $z(0)$ & $\dot{y}(0)$ & $T$ & ($C/B$)-sign \& Floquet multipliers & $\mu_{CZ}$\\
		\hline 3.00014306 &	0.99718487 & 0 & $-$0.00386895 & 4.80 & $(-/+)$ $\varphi_p^3 = 3.217$, $\varphi_s^3 = 0$ & 14$\,\to\,$16 \\
        3.00014210 & 0.99720294 & 0.00036902 & $-$0.00390227 & 4.84 & $(-/+)$ $\varphi_1 = 3.197$, $(-/+)$ $\varphi_2 = 6.282$ & 14 \\
        3.00014182 & 0.99721138 & 0.00041974 & $-$0.00392247 & 4.85 & $(+/+)$ $\lambda = -1.01$, $(-/+)$ $\varphi = 6.281$ & 14 \\
        3.00013324 & 0.99748294 & 0.00114054 & $-$0.00461090 & 5.05 & $(+/+)$ $\lambda = -1.48$, $(-/+)$ $\varphi = 6.089$ & 14 \\
        3.00012084 & 0.99793674 & 0.00156551 & $-$0.00595530 & 5.39 & $(+/-)$ $\varphi_1 = 2.704$, $(-/+)$ $\varphi_2 = 5.684$ & 14 \\
        3.00011446 & 0.99822206 & 0.00164295 & $-$0.00696832 & 5.60 & $(+/-)$ $\varphi_1 = 1.776$, $(-/+)$ $\varphi_2 = 5.499$ & 14 \\
        3.00011146 & 0.99838620 & 0.00163494 & $-$0.00763914 & 5.71 & $0.458 \pm 0.723 i$, $0.624 \pm 0.986 i$ & 14 \\
        3.00010841 & 0.99863569 & 0.00153928 & $-$0.00887793 & 5.82 &	$0.595 \pm 0.212 i$, $1.489 \pm 0.530 i$ & 14 \\
        3.00010824 & 0.99866658 & 0.00151731 & $-$0.00906541 & 5.83 & $(+/+)$ $\lambda_1 = 1.448$, $(-/-)$ $\lambda_2 = 1.786$ & 14 \\
        3.00010818 & 0.99869599 & 0.00149114 & $-$0.00926286 & 5.84 & $(+/+)$ $\lambda_1 = 1.034$, $(-/-)$ $\lambda_2 = 2.566$ & 14 \\
         & & & & & & b-d \\
        3.00010852 & 0.99871851 & 0.00144710 & $-$0.00949590 & 5.84 & $(-/+)$ $\varphi = 6.109$, $(-/-)$ $\lambda = 4.479$ & 13 \\
        3.00010970 & 0.99869031 & 0.00139716 & $-$0.00954957 & 5.83 & $(-/+)$ $\varphi = 6.135$, $(-/-)$ $\lambda = 10.065$ & 13 \\
        3.00011700 & 0.99844373 & 0.00110062 & $-$0.00936425 & 5.75 & $(-/+)$ $\varphi = 6.253$, $(-/-)$ $\lambda = 102.96$ & 13 \\
        3.00012680 & 0.99813343 & 0.00026518 & $-$0.00918878 & 5.65 & $(-/+)$ $\varphi = 6.280$, $(-/-)$ $\lambda = 572.25$ & 13 \\
        3.00012738 & 0.99811615 & 0 & $-$0.00918031 & 5.65 & $\varphi_s^3 = 0$, $(-/-)$ $\lambda_p = 620.1$ & 13$\,\to\,$15
	\end{tabular}}
	\label{green14_saturn_enceladus}
\end{table}
\clearpage
\printbibliography

\end{document}